\documentclass[times,11pt,letter]{article}%
\newif \ifSaveSize \SaveSizetrue
\newif \ifPersonal \Personalfalse% Private version
\usepackage{times,amsmath,amssymb,url}
\newcommand{\ignore}[1]{}
\newcommand{\ensm}[1]{\{ {#1} \}_{N \in \N}~}
\newcommand{\ens}[1]{$\{ {#1} \}_{N \in \N}$~}

\newcommand {\rgs} {random graphs~}

\newcommand {\gnpm} {\mathcal G(N,p)}
\newcommand {\gnp} {$\gnpm$~}
\newcommand {\gnpns} {$\gnpm$}

\newcommand {\kw}{\kwns~}
\newcommand {\kwns}{$k$-wise}
\newcommand {\kwi}{\kwins~}
\newcommand {\kwins}{\kw independent}
\newcommand {\kwice}{\kwicens~}
\newcommand {\kwicens}{\kw independence}
\newcommand {\kwig}{\kwigns~}

\newcommand {\kwigns}{\kw independent graphs}
\newcommand {\kwv}{\kwvns~}
\newcommand {\kwvns}{\kwi variables}

\newcommand {\wig}{\wigns~}
\newcommand {\wigns}{-wise independent graphs}
\newcommand {\wicens}{-wise independence}
\newcommand {\wice}{\wicens~}

\newcommand{\gnkns}{$\gnkm$}

\newcommand{\gnkm}{\mathcal G^{k}(N,p)}
\newcommand{\gnk}{$\gnkm$~}

\newcommand{\gnknm}{\mathcal G^{k(N)}(N,p(N))}
\newcommand{\gnkn}{$\gnknm$~}
\newcommand{\engnk}{\ens{\gnkm} }

\newcommand{\fkt}{\frac k 2}
\newcommand{\fktb} {\fkt} %{\left(\hspace{-0.37ex} \fkt \hspace{-0.2ex}\right)}
\newcommand{\fktbb}{\left(\hspace{-0.37ex} \fkt \hspace{-0.2ex}\right)}
\newcommand{\depth}{\mathit{depth}}
\newcommand{\Edge}{\mbox{\sc edge}}

\newcommand{\ext}{\frac{\log(N)}{\log(1/{p})}}

\newcommand{\kext}{\frac{\log^3 N}{\log^2(1/p)}}

\newcommand {\beq} [1] {\begin {equation} \label{#1}}
\newcommand {\eeq} {\end{equation}}
\newcommand {\bit} {\begin{itemize}}
\newcommand {\eit} {\end{itemize}}
\newcommand {\ben} {\begin{enumerate}}
\newcommand {\een} {\end{enumerate}}
\newtheorem {remark} {Remark}
\newtheorem {Construction} {Construction}
\newcommand \bcons [1] {\begin {Construction} {\bf(#1)}}
\newcommand \econs {\end {Construction}}
\newtheorem {Definition} {Definition}
\newcommand {\bdf} [1] {\begin{Definition} {\bf(#1)}}
\newcommand {\edf} {\end{Definition}}
\newtheorem {Theorem} {Theorem}
\newcommand {\bt} {\begin{Theorem}}
\newcommand {\et} {\end{Theorem}}
\newcommand {\btn} [1] {\begin{Theorem} {\bf(#1)}}
\newcommand {\etn} {\end{Theorem}}
\newtheorem {subTheorem} {Theorem} [Theorem]
\newcommand {\bst} [1] {\begin{subTheorem} {\bf(#1)}}
\newcommand {\est} {\end{subTheorem}}
\newtheorem {Lemma} {Lemma} %[section]
\newcommand {\ble} [1] {\begin{Lemma} {\bf(#1)}}
\newcommand {\ele} {\end{Lemma}}
\newcommand {\blem} {\begin{Lemma}}
\newcommand {\elem} {\end{Lemma}}
\newtheorem {Observation} {Observation}
\newcommand {\bob} [1] {\begin{Observation} {\bf(#1)}}
\newcommand {\eob} {\end{Observation}}
\newtheorem {Claim} {Claim} %[section]
\newcommand {\bcl} {\begin{Claim}}
\newcommand {\ecl} {\end{Claim}}
\newtheorem {subClaim} {Claim} [Claim]
\newcommand {\bscl} [1]{\begin{subClaim}{\bf(#1)}}
\newcommand {\escl} {\end{subClaim}}
\newtheorem {fact} {Fact} %[section]
\newcommand {\bfc} [1] {\begin{fact} {\bf #1}}
\newcommand {\efc} {\end{fact}}
\newtheorem {Corollary} {Corollary}
\newcommand {\bcorl} {\begin{Corollary}}
\newcommand {\ecorl} {\end{Corollary}}
\newtheorem {Condition} {Condition}
\newcommand {\bcond} {\begin{Condition}}
\newcommand {\econd} {\end{Condition}}
\newtheorem {Assumption} {Assumption}
\newcommand {\bas} {\begin{Assumption}}
\newcommand {\eas} {\end{Assumption}}

\newcommand {\mysubsubsection} [1] {\subsubsection{{{#1}}}}%{\subsubsection {{\fbox{#1}}}}
\newcommand {\mysection} [1] {\section{{{#1}}}}%{\section {{\fbox{#1}}}}
\newcommand {\mysubsection} [1] {\subsection{{{#1}}}}%{\subsection {{\fbox{#1}}}}
\newcommand {\myParagraph} [1] {\subsubsection*{#1}}
\newcommand {\SaveSizeParagraph} [1] {\vspace{-1.6mm}\paragraph{#1}}
\newcommand {\SaveaMoreSizeParagraph} [1] {\vspace{-2.7mm}\paragraph{#1}}
\newcommand {\pr} [1] {\Pr \left[ {#1} \right]}

\newcommand {\wrt} {w.r.t.~}
\newcommand {\wip } {with probability~}
\newcommand {\st}{s.t.~}
\newcommand {\as}{a.s.~}
\newcommand {\asns}{a.s.\hspace{-0.2ex}}
\newcommand {\etal} {et al.~\hspace{-0.1ex}}
\newcommand {\qed} {$\blacksquare$}
\newcommand {\qedm} {~\blacksquare}
\newcommand {\eqdef} {{~ \stackrel{\mathrm{def}} {=} ~}}
\newcommand {\poly} {\mathit{poly}}
\newcommand {\Z} {\mathbb{Z}}
\newcommand {\N} {\mathbb{N}}
\newcommand {\lb} {\left(}
\newcommand {\rb} {\right)}
\newcommand {\lsb} {\left[}
\newcommand {\rsb} {\right]}
\newcommand {\ceil}[1] {\lceil {#1} \rceil}
\newcommand {\floor}[1] {\lfloor {#1} \rfloor}
\newcommand {\into} {\rightarrow} %function F:A \into B
\newcommand {\eps} {\epsilon}
\newcommand {\opmo} {{\scriptstyle (1 \pm o(1))}}
\newcommand {\opo} {{\scriptstyle (1 + o(1))}}
\newcommand {\omo} {{\scriptstyle (1 - o(1))}}
\newif \ifCiteName \CiteNamefalse %cite [GGM] or [3]
\newcommand {\myBib} [2]{\ifCiteName \bibitem[#1]{#2}\else \bibitem{#2}\fi}
\newcommand \myTitle {\emph}
\newcommand \myCite {}
\newcommand \andCoAuthersBib {, }
\newcommand \andCoAuthers {and }
\newcommand {\ER} {Erd\"{o}s \andCoAuthers R\'{e}nyi~}

\newcommand {\Erd} {Erd\"{o}s~}

\newcommand {\Bol} {Bollob\'{a}s~}
\begin{document}
\thispagestyle{empty}
\title{\bf{k-wise independent random graphs}}
\author{Noga Alon\thanks{
Schools of Mathematics and Computer Science,
Sackler Faculty of Exact Sciences, Tel Aviv
University, Tel Aviv 69978,
Israel, and IAS, Princeton, NJ 08540, USA.
Email:~{\texttt{nogaa@post.tau.ac.il.}}
Research supported in part by the Israel Science
Foundation and by a USA-Israeli BSF
grant.}
\and
Asaf Nussboim\thanks{Department of Computer Science and
Applied Mathematics, Weizmann Institute of Science,
Rehovot, Israel.
Email:~{\texttt{asaf.nussbaum@weizmann.ac.il.}}
{Partly supported by a grant from the Israel Science
Foundation.}}}
\date{}\maketitle
%=====================================================
%                      Abstract
%=====================================================
\begin {abstract}
We study the \kwi relaxation of the usual model \gnp of random
graphs where, as in this model,
$N$ labeled vertices are fixed and each edge
is drawn with probability $p$,
however, it is only required that the distribution of
any subset of $k$ edges is independent.
This relaxation can be relevant in modeling phenomena
where only \kwice is assumed to hold, and is also useful
when the relevant graphs are so huge
that handling \gnp graphs becomes infeasible, and
cheaper random-looking distributions (such as \kwi
ones) must be used instead. Unfortunately, many well-known
properties of random graphs in \gnp are global, and it is
thus not clear if they are
guaranteed to hold
in the \kwi case.
We explore the properties of \kwig by providing
upper-bounds and lower-bounds on the amount of
independence, $k$, required for maintaining the main
properties of \gnp graphs: connectivity,
Hamiltonicity, the connectivity-number, clique-number
and chromatic-number and the appearance of fixed
subgraphs.
Most of these properties are shown to be captured by
either constant $k$ or by some $k=\poly(\log(N))$ for
a wide range of values of $p$, implying that random looking graphs
on $N$ vertices can be generated by a seed of size
$\poly(\log(N))$. The proofs combine combinatorial, probabilistic
and spectral techniques.

\end {abstract}
\thispagestyle{empty}
\newpage \setcounter{page}{1}
%=======================
\mysection{Introduction}
%=======================
We study the \kwi relaxation of the usual model \gnp of random graphs
where, as in this model, $N$ labeled vertices are fixed and each
edge is drawn with probability (w.p., for short)~$p=p(N)$, however, it is only
required that the distribution of any subset of $k$
edges is independent (in \gnp all edges are
mutually independent).
These \kwig are natural combinatorial objects that
may prove to be useful in modeling scientific phenomena
where only \kwice is assumed to hold. Moreover,
they can be used when the relevant graphs are so huge, that handling \gnp graphs is infeasible, and cheaper random-looking distributions must be used instead.
However, what happens when the application that uses
these graphs (or the analysis conducted on them)
critically relies on the fact that \rgs are, say,
almost surely connected? After all, \kwice is defined via `local' conditions, so isn't it possible that
\kwig will fail to meet `global' qualities like
connectivity? This motivates studying which global attributes
of \rgs are captured by their \kwi counterparts.

Before elaborating on properties of \kwig we provide
some background on \kwicens, on properties of random
graphs, and on the emulation of huge random graphs.

\mysubsection
%================================
{Emulating Huge Random Graphs}
%================================
Suppose that one wishes to conduct some simulations on
random graphs. Utilizing \gnp graphs requires
resources polynomial in $N$, which is infeasible
when $N$ is huge (for example, exponential in the
input length, $n$, of the relevant algorithms).
A plausible solution is to replace \gnp by a
cheaper `random looking' distribution $\mathcal G_N$.
To this end, each graph $G$ in the support of
$\mathcal G_N$ is represented by a very short binary
string (called seed) $s(G)$, \st evaluating edge
queries on $G$ can be done efficiently when $s(G)$ is
known; Then, sampling a graph from $\mathcal G_N$ is
done by picking the seed uniformly at random.

Goldreich \etal were the first to address this scenario
in \cite{ggn}. They studied emulation by computationally
pseudorandom graphs, that are indistinguishable from \gnp from the view of
any $\poly(\log(N))$-time algorithm that inspects
graphs via edge-queries of its choice. They considered
several prominent properties of
\gnp graphs, and constructed computationally
pseudorandom graphs that preserve many of those
properties (see the final paragraph of Section \ref{compareGGN}).

We consider replacing \rgs by \kwi ones. The latter
can be sampled and accessed using only
$\poly(k\log(N))$-bounded resources. This is achieved
thanks to efficient constructions of discrete \kwi
variables by Joffe \cite{jof}, see also Alon, Babai and
Itai \cite{abi}:
the appearance of any potential edge in the graph is
simply decided by a single random bit (that has
probability $p$ to attain the value 1).
Such \kwig were used by Naor \etal \cite{nnt} to
efficiently capture arbitrary first-order properties
of huge \gnp graphs (see Section \ref{FO_paper}).

\mysubsection
%================================================================
{${\bf k}$-Wise Independent Random Variables}
%================================================================
Distributions of discrete \kwi variables play an
important role in computer science. Such distributions
are mainly used for de-randomizing algorithms (and for
some cryptographic applications). In addition, the
complexity of constructing \kwi variables was studied
in depth, and in particular, the aforementioned constructions
\cite{jof,abi} (based on degree $k$ polynomials over finite fields)
are known to provide essentially the smallest possible sample spaces.
Our work is, however, the first systematic
study of {\em combinatorial properties} of \kwi
objects. Properties of various other \kwi objects
(mainly percolation on $\Z^d$ and on Galton-Watson
trees) were subsequently explored by Benjamini,
Gurel-Gurevich and Peled \cite{bggp}.

\mysubsection
%================================================================
{The Combinatorial Structure of Random Graphs}
%================================================================
%\ifSaveSize
What are the principal attributes of \rgs that \kwi ones should maintain? Most theorems that manifest the remarkable structure of \rgs state that certain properties occur either almost surely (\as for
short), or alternatively hardly ever, (namely, \wip
tending either to 1 or to 0 as $N$ grows to $\infty$).
These results typically fall into one of the following
categories.

%===========================================================
\SaveSizeParagraph{Tight concentration of measure.}
%===========================================================
A variety of prominent random variables (regarding
random graphs) \as attain only values that are {\em
extremely close} to their expectation. For instance,
random graphs (with, say, constant $p$) \as have
connectivity number $\kappa=\opmo pN$, clique number
$c=\opmo \frac {2\log(pN)}{\log(1/p)}$ (\Bol
\andCoAuthers \Erd\cite{be}, Matula \cite{mat}, Frieze
\cite{fri}) and chromatic number $\chi=\opmo \frac
{N\log(1/1-p)}{2\log(pN)}$ (\Bol \cite{bolChi},
{\L}uczak \cite{l}).

%===========================================================
\SaveSizeParagraph{Thresholds for monotone
properties.}
%===============================================================
For a given monotone increasing\footnote{Namely, any
property closed under graph isomorphism and under
addition of edges.} graph property $T$, how large
should $p(N)$ be for the property to hold \asns?
This question had been settled for many prominent properties such as connectivity (\ER
\cite{erConn}), containing a perfect matching (\ER
\cite{erMatch1,erMatch2,erMatch3}), Hamiltonicity
(P\'{o}sa \cite{pos}, Kor\v sunov \cite{kor},
Koml\'{o}s \andCoAuthers Szemer\'{e}di \cite{ks0}), and
the property of containing copies of some fixed
graph $H$ (\ER \cite{erGiant}, \Bol \cite{bolsub}).
For these (and other) graph properties the sufficient
density (for obtaining the property) is surprisingly
small, and moreover, a threshold phenomenon occurs when
by `slightly' increasing the density from
$\underline{p}(N)$ to $\overline{p}(N)$, the
probability that $T$ holds dramatically changes from
$o(1)$ to $1-o(1)$.%
\footnote{Thresholds for prominent properties are
often so sharp that $\overline{p}=(1+o(1))
\underline{p}$. Somewhat coarser thresholds were (later)
established for {\em arbitrary} monotone properties by
\Bol \andCoAuthers Thomason \cite{bt}, and by Friedgut
\andCoAuthers Kalai \cite{fk}.}
Thus, good emulation requires the property $T$ to be guaranteed at densities as close as possible to the true \gnp threshold.

%===========================================================
\SaveSizeParagraph{Zero-one laws.}
%===========================================================
These well known theorems reveal that {\em any} first-order
property holds either \as or hardly ever for \gnpns. A
first-order property is any graph property that can be
expressed by a single formula in the canonical
language where variables stand for vertices and the
only relations are equality and adjacency (e.g.~``having an isolated vertex" is specified
by $\exists x \forall y \neg \Edge (x,y)$).
These Zero-one laws hold for any fixed $p$ (Fagin \cite{fag}, Glebskii, Kogan, Liagonkii \andCoAuthers Talanov \cite{gklt}), and whenever $p(N)=N^{-\alpha}$ for a fixed irrational $\alpha$ (Shelah \andCoAuthers
Spencer \cite {ss}).

\mysection
%================================================================
{Our Contribution} %================================================================

We investigate the properties of \kwig by providing
upper bounds and lower bounds on the `minimal' amount of
independence, $k_T$, required for maintaining the main properties $T$
of random graphs.\ifPersonal\footnote{When $k={{N \choose 2}}$ we get the original \gnp graphs, while for, say, $k=1$ (e.g., the case where \wip $p$ the graph is complete and otherwise it is empty) the resulting graphs have very little in common with \gnp graphs.} \else{ }\fi
The properties considered are:
connectivity, perfect matchings, Hamiltonicity, the
connectivity-number, clique-number and
chromatic-number and the appearance of copies of a
fixed subgraph $H$.
We mainly establish upper bounds on $k_T$ (where arbitrary \kwig are shown to exhibit the property $T$) but also lower bounds (that provide specific constructions of \kwig that fail to preserve $T$).
Our precise results per each of these properties are
discussed in Section \ref{intro_list_results},
and proved in Section \ref{claims_&_proofs} (and the Appendices).
Interestingly, our results reveal a deep difference
between \kwice and almost \kwice (a$.$k$.$a$.$
$(k,\eps)$-\wicens\footnote{$(k,\eps)$-\wice means
that the joint distribution of any $k$ potential edges
is only required to be within small statistical
distance $\epsilon$ from the corresponding
distribution in the \gnp case.}).
All aforementioned graph properties are guaranteed by
\kwice (even for small $k=\poly(\log(N))$), but are
strongly violated by some almost \kwig - even when
$k=N^{\Omega(1)}$ is huge and $\eps=N^{-\Omega(1)}$ is
tiny.
For some properties of random graphs, $T$, our results
demonstrate for the first time how to efficiently
construct random-looking distributions on huge graphs
that satisfy $T$.

\SaveSizeParagraph
%============================================================
{Our Techniques \& Relations to Combinatorial Pseudorandomness.}
%============================================================
\label{intro_jumbled}

For positive results (upper bounding $k_T$), we note that the original proofs that establish properties of \gnp
graphs often fail for \kwigns. These proofs use a union
bound over $M=2^{\Theta(N)}$ undesired events, by
giving a $2^{-\Omega(N)}$ upper-bound on the
probability of each of these events.%
\footnote{For instance \wrt connectivity, $M$ is the number of choices for partitioning the vertices into 2 disconnected components.%
\ifPersonal{~For perfect matchings, $M$ counts the
number of subsets of vertices that defy a sufficient
condition for the matching (by either Hall's Theorem or
Tutte's Theorem).
For the chromatic number, $M$ is the (a-priory) number of possible sub-graphs from which a greedy coloring
algorithm might choose a large independent set that is
colored with a single color.}\fi}
Unfortunately, there exist $\poly(\log(N))$-\wig
\ifPersonal that can be constructed consuming only
$\Theta(k \log(N))$ random bits. Hence, any event that
occurs with positive probability (in these
constructions), must hold \wip $\geq 2^{-\Theta(k
\log(N))}$ which is $\gg \frac 1 M$ when $k \leq
\poly(\log(N))$. \else where any event that occurs
with positive probability, has probability $\geq
2^{-o(N)}$. \fi
Therefore, directly `de-randomizing' the original proof fails, and alternative arguments (suitable for the $k$-wise independent case) are provided.

In particular, many properties are inferred via a variant of
Thomason's notion of `jumbledness' \cite{t}
(mostly known in its weaker form as quasirandomness or
pseudorandomness, as defined by Chung, Graham \andCoAuthers Wilson
\cite{cgw}, and related to the so called Expander Mixing Lemma and the
pseudo-random properties of graphs
that follow from their spectral properties, see \cite{ac}).
For our purposes, $\alpha$-jumbledness means that (as expected in
\gnp graphs) for all vertex-sets $U,V$, the number of
edges that pass from $U$ to $V$ should be $p|U||V|
\pm \alpha \sqrt{|U||V|}$.
Jumbledness and quasirandomness had been studied extensively (see
\cite{ks} and its many references), and serve in Graph Theory as
{\em the} common notion of resemblance to random graphs. In
particular, \gnp graphs are known to exhibit (the best
possible) jumbledness parameter,
$\alpha=\Theta(\sqrt{pN})$.
One of our main results (Theorem
\ref{kwise_gives_jumbl}) demonstrates that \kwice for
$k=\Theta( \log(N))$ is stronger than jumbledness, in
the sense that it guarantees the optimal
$\alpha=\Theta(\sqrt{pN})$ even for tiny densities
$p=\Theta(\frac{\ln(N)}N)$.
Therefore, prominent properties of \kwig can be directly deduced from properties of jumbled graphs.

Proving Theorem \ref{kwise_gives_jumbl} exploits a
known connection between jumbledness and the
eigenvalues of (a shifted variant of) the adjacency
matrix of graphs, following the approach in \cite{ac}.
In particular, the analysis of Vu (\cite{vu},
extending \cite{fk2}) regarding the eigenvalues of
random graphs is strengthened, in order to achieve
optimal eigenvalues even for smaller densities $p$
than those captured by \cite{vu}. This improvement
implies, among other results, the remarkable fact that
\kwig for $k=\Theta(\log(N))$ preserve (up to constant
factors) the \gnp sufficient density for connectivity.

\SaveSizeParagraph
%============================================================
{More on Techniques \& Relations to Almost $k$-Wise Independence.} %============================================================
\label{almost_k_wise}

For negative results (producing random-looking graphs
that defy a given property $T$ of random graphs), the
\cite{ggn,my_thesis} approach is to first construct
some random-looking graph $G$, and later to `mildly'
modify $G$ \st $T$ is defied. This is done \wrt all graph properties considered here. For instance, the
modification of choosing a random vertex and then deleting all it's edges violates connectivity while preserving computational pseudorandomness.
Unfortunately, such modifications \ifSaveSize \else are useless in our context because they \fi fail to preserve
\kwice \ifPersonal and it is typically not clear how
to modify the graph for the second time to regain \kwice \else \fi (the resulting graphs are only almost \kw independent).
In contrast, most of our negative results exploit the fact
that some constructions of \kwi bits produce
strings with significantly larger probability than in
the completely independent case. This is translated
(by the construction in Lemma
\ref{kwise_construction_unexpected_sub_graphs}) to the
unexpected appearance of some subgraphs (in \kwigns):
either huge independent-sets inside dense graphs or
fixed subgraphs inside sparse graphs.

\SaveSizeParagraph
%============================================================
{Comparison with Computational Pseudorandomness.} %============================================================
\label{compareGGN}

Finally, \kwice guarantees all random graphs'
properties that were met by the (specific)
computationally pseudorandom graphs of
\cite{ggn,my_thesis}\ifSaveSize. \else \footnote{A
single exception is preserving the precise \gnp
upper-bound on the chromatic number only up to a
constant factor, whereas no such factors are
introduced by \cite{ggn,my_thesis}.}. \fi
\ifSaveSize In addition, only \kwice captures \else
The list of random graph properties captured by \kwice
but not by \cite{ggn,my_thesis} includes \fi (i)
arbitrary first-order properties of \gnp graphs, (ii)
high connectivity, (iii) strongest possible parameters
of jumbledness, and (iv) almost regular $(1\pm
o(1))pN$ degree for all vertices, and $(1\pm
o(1))p^2N$ co-degrees for all vertex pairs.
Importantly, all this holds for any \kwigns,  (and in
particular for the very simple and efficiently
constructable ones derived from \cite{jof,abi}),
whereas the \cite{ggn,my_thesis}'s approach requires
non-trivial modifications of the construction per each
new property.

\mysection
%================================================================
{Combinatorial Properties of ${\bf k}$-Wise Independent %================================================================
Graphs} \label{intro_list_results}

%\ifSaveSize
We now survey our main results per each of the
aforementioned graph properties $T$. Typically our
arguments establish the following tradeoff: the
smaller $p$ is, the larger $k$ should be to maintain
$T$.
Given this tradeoff we highlight minimizing $k$ or, alternatively, minimizing $p$. The latter is
motivated by the fact that the \gnp threshold for
many central properties occurs at some $p^*\ll 1$.
Minimizing $p$ is subject to some reasonable choice of
$k$, which is $k\leq \poly(\log(N))$. Indeed, as the
complexity of implementing \kwig is $\poly(k
\log(N))$, we get efficient implementations whenever
$k \leq \poly(\log(N))$ even when the graphs are huge
and $N=2^{\poly(n)}$.
\footnote{Accessing the graphs via edge-queries is
adequate only when $p \geq n^{-\Theta(1)}$ - otherwise
\as no edges are detected by the $\poly(n)$ inspecting
algorithm. For smaller densities our study has thus mostly a
combinatorial flavor.}

\ifSaveSize{\vspace{-2.7mm}}\fi \mysubsection
%==================================================
{Connectivity, Hamiltonicity and Perfect Matchings
%==================================================
(see Section \ref{section_conn_Ham})}

The well known sufficient \gnp density for all these
properties is $\sim \frac {\ln(N)}N$. For
connectivity, this sufficient density is captured (up
to constant factors) by all $\log(N)$-\wigns. Even
$k=4$ suffices for larger densities $p \gg N^{-\frac 1
{2}}$.
Based on Hefetz, Krivelevich \andCoAuthers Szabo's
\cite{hks}, Hamiltonicity (and hence perfect
matchings) are guaranteed at $p\geq \frac
{\log^2(N)}{N}$ with $k\geq 4\log(N)$, and at $p \geq
N^{-\frac 1 {2}+o(1)}$ with $k\geq 4$.
On the other hand, some pair\wig are provided that
despite having constant density, are still \as
disconnected and fail to contain any perfect matching.
\ifPersonal Note that by changing $k$ to 4, and the
density to $p=1/2 + N^{-\Theta(1)}$, all properties are
ensured again (since then all the degrees are $\leq
\frac N 2$). More importantly, note that the combination
of these results immediately rule out the possibility of
establishing a general threshold phenomenon (for all
monotone properties) when $k=2$.} \fi

\ifSaveSize{\vspace{-2.7mm}}\fi \mysubsection
%=================================
{High Connectivity
%=================================
(see Section \ref{section_high_conn})}

The connectivity number, $\kappa(G)$, is the largest
integer, $\ell$, \st any pair of vertices is connected
in $G$ by at least $\ell$ internally vertex-disjoint
paths. Since a typical degree in a random graph is
$(1\pm o(1)) pN$, it is remarkable that \gnp graphs
achieve $\kappa = (1\pm o(1)) pN$ \asns.
Surprisingly, such optimal connectivity is guaranteed
by $\Theta(\log(N))$\wice whenever $p \geq
\Theta(\frac {\log(N)}{N})$.

\ifSaveSize{\vspace{-2.7mm}}\fi \mysubsection
%=======================================================
{Cliques and Independent-Sets
%=======================================================
(see Appendix \ref{appendix_indp_num})}

For $N^{-o(1)} \leq p \leq 1-N^{o(1)}$ the
independence number, $I$, of random graphs has \as
only two possible values: either $S^*$ or $S^*+1$ for
some $S^*\sim \frac {2\log(pN)}{\log(1/(1-p))}$.
This remarkable phenomenon is observed to hold by
$\Theta(\log^2(N))$-\wice whenever $p$ is
bounded away from 0.
On the other hand, \kwig are provided with $k= \Theta
\lb \frac {\log(N)}{\log \log(N)}\rb$ where $I \geq
(S^*)^{1+\Omega(1)}$ \as (for $k=\Theta(1)$, even huge
$N^{\Omega(1)}$ independent-sets may appear).
For smaller densities, random graphs \as have $I \leq
O(p^{-1} \log(N))$, while $\Theta(\log(N))$\wice gives
a weaker, yet useful, $I \leq O(\sqrt{\frac N p})$
bound whenever $p\geq \Omega (\frac{\log(N)}N)$.
By symmetry (replacing $p$ with $1-p$), analogous
results to all the above hold for the clique number as
well.
Discussing the clique- and independence-number is
deferred to the appendices since the main
relevant techniques here are demonstrated elsewhere in the paper.

\ifSaveSize{\vspace{-2.7mm}}\fi \mysubsection
%=======================================================
{Coloring
%=======================================================
(see Section \ref{section_color})}

For $1/N \ll p \leq 1-\Omega(1)$, the chromatic number $\chi$ of random graphs
is a.s. $(1+o(1))\frac
{N\log(1/1-p)} {2\log(pN)}$.
This \gnp lower-bound on $\chi$ is observed to hold
for any $(\log(N))^{\Theta(1)}$\wig with moderately
small densities $p\geq (\log(N))^{-\Theta(1)}$.
More surprisingly, $k=\Theta(\log(N))$ suffices to
capture a similar upper-bound (even for tiny densities
$p=c \log(N)/N$).
%$p=N^{-\Theta(1)}$).
%
This upper-bound is based on Alon, Krivelevich
\andCoAuthers Sudakov's \cite{aks}, \cite{aks1} and on Johansson's
\cite{joh}.

\ifSaveSize{\vspace{-2.7mm}}\fi \mysubsection
%=====================================================================
{Thresholds for the Appearance of Subgraphs
%=====================================================================
(see Section \ref{section_sub_graphs})} \label{intro_H_copies}

For a fixed (non-empty) graph $H$, consider the
appearance of $H$-copies ({\em not necessarily} as an
induced subgraph) in either a random or a \kwi
graph.\ifPersonal\footnote{
When $H$ is empty, the question is trivial because every
graph $G$ (of sufficient order) contains $H$ copies.
This formally translates to $\rho=\infty$ and to
$p^*=0$.}\else~\fi
The \gnp threshold for the occurrence of $H$ sub-graphs
lies at $p^*_H \eqdef N^{-\rho}$, where the constant
$\rho=\rho(H)$ is the minimum, taken over all subgraphs
$H'$ of $H$ (including $H$ itself), of the ratio $\frac
{v(H')} {e(H')}$ (here, $v(H')$ and $e(H')$ respectively
denote the number of vertices and edges in $H'$).
Thus, no $H$-copies are found when ${p}\ll p^*$, while for any ${p}\gg p^*$, copies of $H$ abound (\ER \cite{erGiant}, \Bol \cite{bolsub}).
\ifPersonal %-----------------------------------------
\begin{remark} The ratio $\rho$ decides the
threshold, since it is minimized for the most `dense'
subgraph $H'$. Intuitively, this $H'$ is the sub-graph
less likely to appear in a \gnp graph, so $H$-copies
\as appear iff $H'$-copies \as appear.\end{remark}

\begin{remark} Note that $\frac 1 {v-1} \leq \rho \leq
2$, where the lower bound is tight (only) for cliques,
and the upper bound is tight iff $H$ contains only
disjoint edges. As $v=v(H)$ grows, most graphs have
$e(H)=\Theta(v^2)$, so $\rho=\Theta(1/v)$.\end{remark}

\fi
%\ifPersonal Remarks regarding \rho
%-----------------------------------------------------
%
For any graph $H$, this \gnp threshold holds whenever
$k\geq \Theta (v^4(H))$, but as $k$ is decreased to
$\floor{\frac 2 {\rho}}$, the \gnp threshold is defied:
much sparser graphs exist where $p \ll p^*_H$ and yet
copies of $H$ are \as found.
In particular, when $e(H)\geq \Omega(v^2(H))$,
\ifPersonal then $\rho \leq O(1/v(H))$, so \fi the
threshold violation occurs at $k= \Omega(v(H))$.
\ifPersonal Note that our negative results are optimal
in the sense that they capture all possible graphs $H$
(see footmark inside statement of Theorem
\ref{sub_graphs_defy}). The optimality of $k$ is
discussed in remark \ref{remark_rho_cond_reasonable}.
\fi

\ifSaveSize{\vspace{-2.7mm}}\fi \mysubsection
%================================================================
{First Order Zero-One Laws (Previous Results)}
%================================================================
\label{FO_paper}

\ifSaveSize Naor \etal \cite{nnt} have recently \else a recent study (joint work of the second author with Naor \andCoAuthers Tromer \cite {nnt}) \fi studied capturing arbitrary depth-$D(N)$ properties.
These are graph properties expressible by a sequence of first-order formulas $\Phi=$ $\{\phi_N\}_{N \in \mathbb N}$, with quantifier depth $\depth(\phi_N) \leq D(N)$%
\ifSaveSize. \else
(e.g.~``having a clique of size $t(N)$" can be specified
by $\phi_N = \exists x_1 ... \exists x_{t(N)}
\bigwedge_{i \neq j} ((x_i \neq x_j) \bigwedge
\Edge(x_i,x_j))$). \fi
A `threshold' depth function $D^*\sim\ext$ was
identified \st arbitrary \kwig resemble \gnp graphs \wrt all
depth $D^*$ properties. The underlying resemblance-definition is in
fact so strong, that even \gnp graphs cannot achieve resemblance
to themselves \wrt properties of higher depth. On the other hand,
\kwig were shown to defy some \gnp properties of depth
$\Theta(\sqrt{k\log(N)}+\log(N))$.
These results are incomparable to the ones in the current paper,
since most of the graph properties studied here require larger
depth than $D^*$.
\ifPersonal\footnote{ Note that the naive recursive
formula for connectivity has depth $3 \log(N)$ but
exponential size, while Savitch's Theorem \cite{sav}
gives linear size and $5 \log(N)$ depth. On the other
hand, by a simple Ehrenfeucht-game argument, depth $2
\log(N)-\Theta(1)$ is insufficient for connectivity.
For optimal Ramsey graphs depth $2 \log(N)$ suffices,
which is again tight by some trivial Ehrenfeucht-game
argument.}
Thus, a more direct approach is needed to establish such
properties for \kwigns.\footnote{Note that, in addition,
direct proofs may reduce the amount of independence $k$
required. For instance for appearance of subgraphs of
order $v$, then $k=\Theta(v^4)$ suffices instead of $k
\sim \kext$.}\fi

\mysection{Preliminaries}%%\SaveaMoreSizeParagraph
%============================================================
\paragraph{Asymptotics.}
%============================================================
Invariably, $k: \N \into \N$, while
$p,\eps,\delta,\gamma,\Delta: \N \into (0,1)$. We
often use $k,p,\eps,\delta,\gamma,\Delta$ instead of
$k(N),p(N),\eps(N),\delta(N),\gamma(N),\Delta(N)$.
Asymptotics are taken as $N \rightarrow \infty$, and
some inequalities hold only for sufficiently large
$N$. The $\floor{\cdot}$ and $\ceil{\cdot}$ operators
are ignored whenever insignificant for the asymptotic
results. Constants $c,\bar c$ are not optimized in
expressions of the form $k=c\log(N)$ or
$p=(\log(N))^{\bar c}/N^{\Delta}$, whereas the
constant $\Delta$ is typically optimized.

%===========================================================
\SaveaMoreSizeParagraph{Subgraphs.}
%============================================================
For a graph $H$, let $v(H)$ and $e(H)$ denote the number
of vertices and edges in $H$. For vertex sets $U,V$ let
$e(U,V)$ denote the number of edges that pass from $U$ to
$V$ (if $S=U \bigcap V \neq \emptyset$, then any
internal edge of $S$ is counted twice). Similarly, we
let $e(U)=e(U,U)$.

%===========================================================
\SaveaMoreSizeParagraph{Random and ${\bf k}$-Wise
Independent Graphs.}
%============================================================
Throughout, graphs are simple, labeled and undirected.
Given $N,k,p$ as above then \gnkn (or \gnk for short)
denotes some distribution over the set of graphs with
vertex set $\{1,...,N\}$, where each edge appears
w.p.~$p(N)$, and the random variables that indicate the
appearance of any $k(N)$ potential edges are mutually
independent. We use the term `\kwigns' for a sequence of
distributions \ens{\gnkm} indexed by $N$.

%===========================================================
\SaveaMoreSizeParagraph{Almost Sure Graph Properties.}
%============================================================
A graph property $T$, is any property closed under
graph isomorphism. We say that `$T$ holds a.s.~(almost
surely) for \gnkns' or that (abused notation) `$T$
holds for \gnkns' whenever $\Pr_{\gnkm}[T]$ {}
$\stackrel {N \into \infty} {\longrightarrow} 1$. Similar terminology is used for \gnp graphs.

%===========================================================
\SaveaMoreSizeParagraph{Monotonicity in $(\bf{k,p})$.}
%============================================================
Since $\bar k$-\wice implies $k$-\wice for all $\bar k >
k$ we may state claims for arbitrary $k\geq k'$ but
prove them only for $k=k'$.
When establishing monotone increasing properties we
often state claims for arbitrary $p\geq p'$ but prove
them only for $p=p'$.
The latter is valid since for any $N,k,p>p'$, the
process of sampling from any (independent) $\mathcal G^k(N,p)$,
$\mathcal G^k(N,p'/p)$ distributions and defining the
final graph with edge-set being the intersection of the
edge-sets of the two sampled graphs, clearly results in a $\mathcal G^k(N,p')$ distribution.
\ifPersonal Assuming that $k=k',p=p'$ is used only to
ensure the $\frac {M-k}k \mu(1-\mu) \geq 1$ condition in
Lemma \ref{kwise_Chebyshev_bound} (the $k \leq \frac N
3$ condition holds anyway, whenever the Lemma is
meaningful).
Finally, lower-bounds on $p$ are often provided only for
sake of clarity (these lower-bounds are redundant as
they follow immediately from the companion bounds on
$\eps$ or $\gamma$).\fi

%===========================================================
\SaveaMoreSizeParagraph{${\bf k}$-Wise Independent
Random Variables.}
%============================================================
The term `$(M,k,p)$-variables' stands for any $M$
binary variables that are $k$-wise independent with
each variable having probability $p$ of attaining
value 1. Lemma \ref{modify_JoffeCG_lem} (proved in
Section \ref{modify_JoffeCG_sec}) adjusts the known
construction of discrete \kwv of \cite{jof},\cite{cg}, \cite{abi} to
provide $(M,k,p)$-variables that induce some
predetermined values with relatively high
probability. Throughout, $e_1$ and $e_0$ resp.~denote
the number of edges and non-edges in a graph $H$.

\blem \label{modify_JoffeCG_lem} Given $0<p<1$ with
binary representation $p=0.b_1...b_{\ell}$, and
natural numbers $e_0,e_1,M$ satisfying $e_0+e_1\leq
M$, let $F= \max \{2^{\ceil{\log_2M}},2^{\ell}\}$.
Then there exists $(M,k,p)$-variables \st
$\Pr[A]=F^{-k}$, where $A$ denotes the event that the
first $e_0$ variables receive value 0 while the next
$e_1$ variables receive value 1. \elem

%===========================================================
\SaveaMoreSizeParagraph{Tail Bounds for ${\bf k}$-Wise
Independent Random Variables.}
%============================================================
The following strengthened version of standard tail
bounds (proved in Section \ref{sec_inequal})
translates into smaller densities $p$ for which
monotone graph properties are established for \kwigns.

\blem \label{kwise_Chebyshev_bound} Let $X=\sum_{j=1}^M
X_{j}$ be the sum of \kwi binary variables where
$\Pr[X_j=1]=\mu$ holds for all $j$. Let $\delta > 0$,
and let $k$ be even \st $\frac {M-k}k \mu(1-\mu) \geq
1$. Then
$$\Pr[|X-\mathbb E(X)| \geq \delta \mathbb E(X)] \leq
\left [\frac {2k(1-\mu)} {\delta^2 \mu M}\right] ^{\fktb}.$$
\elem

%======================================================

\mysection{The properties of ${\bf k}$-wise
independent graphs} \label{claims_&_proofs}
\mysubsection{Degrees, Co-Degrees and Jumbledness}
\ifPersonal For future purposes, we first establish some
fundamental random graphs' properties, such as having
almost regular degrees and co-degrees as well as
achieving strong jumbledness.\fi

%========================================================
%                       Degrees
%========================================================

\ble{Achieving almost regular degrees}
\label{kwise_gives_deg}
In all \kwig \engnk it \as holds that all vertices have degree~~$p(N-1)(1 \pm \eps)$ whenever $N \big[\frac {3k} {\eps^2
pN}\big]^{\floor{k/2}} \longrightarrow 0,$ and in
particular when either \ben
\item
$k\geq 4$,~$N^{-1/2} \ll p \leq 1-\frac 5 N$, and $1 \geq \eps \gg p^{-1/2}N^{-1/4};$~~or
\item
$k\geq 4\log(N)$,~$\frac{25\log(N)}{N} \leq p \leq
1-\frac {5 \log(N)} N$, and $1 \geq \eps \geq
\sqrt{\frac{25\log(N)}{p N}}.$ \een \ele
%
%----------------------------------------------------
%
\noindent{\bf Proof.} Fix a vertex $v$, and let $X_w$ be
the random variable that indicates the appearance of the
edge $\{v,w\}$ in the graph. Thus, the degree of $v$ is
$X=\sum_{w\neq v}X_w$. Since $X$ is the sum of
$(N-1,k,p)$-variables, Lemma \ref{kwise_Chebyshev_bound}
implies that the probability that $v$ has an unexpected
degree $X \neq p(N-1)(1\pm\eps)$ is bounded by
$\big[\frac {3k} {\eps^2 pN}\big]^{\floor{k/2}}.$
Applying a union-bound over the $N$ possible vertices
$v$, gives that the probability of having {\em some}
vertex with unexpected degree is bounded by $N
\big[\frac {3k} {\eps^2 pN}\big]^{\floor{k/2}},$ which
vanishes for the parameters in items 1 and 2. \qed

%========================================================
%                       Co-Degrees
%========================================================

\ble{Achieving almost regular
co-degrees} \label{kwise_gives_codeg}
In all \kwig \engnk it \as holds that all vertex pairs
have co-degree~~$p^2(N-2)(1\pm \gamma)$ whenever
either \ben
\item
$k\geq 12$,~$N^{-\frac 1 6} \ll p \leq 1-\frac {13} N$, and $1 \geq \gamma \gg p^{-1} N^{-\frac 1 6};$~~or
\item
$k\geq 12\log(N)$,~$\sqrt{\frac{73 \log(N)}{N}} \leq p
\leq 1-\frac {13 \log(N)} N$~~and~~$1 \geq \gamma \geq
\sqrt{\frac {73 \log(N)} {p^2N}}.$ \een \ele
%
%Complete Proof ---------------------------------------
%
\ifPersonal{%
\noindent{\bf Proof.} We prove item 1 for $k=12$, item
2 for $k=12\log(N)$ and item 3 for $k=4$. The claim
follows for $\bar k > k$, since $\bar k$-\wice implies
$k$-\wicens.
Fix a vertex pair $\{u,v\}$, and let $X_w$ be the random variable indicating the appearance of both edges $\{u,w\}$, $\{v,w\}$ in the graph. Thus the co-degree of $\{u,v\}$ is $X=\sum_{w\neq u,v}X_w$. %
Now $X$ is the sum of $N-2$ binary variables that are
$\lb k/2 \rb$-wise independent and each has probability
of success $p^2$. Thus, Lemma
\ref{kwise_Chebyshev_bound} combined with a union-bound
over the $\binom N 2 $ possible vertex pairs $\{u,v\}$,
ensures that the probability that there exist {\em some}
vertex-pair with unexpected co-degree $X \neq
p^2(N-2)(1\pm\gamma)$ is bounded by $\binom N 2 \left
[\frac {3k} {\gamma^2 p^2N}\right]^{\frac k 4}.$
Clearly, the latter expression vanishes for our choice
of parameters. \qed
%
%Proof Sketch-----------------------------------------
%
\else \noindent{\bf Proof.} The proof is completely
analogous to that of Lemma \ref{kwise_gives_deg}. Here the
union-bound is over all $\binom N 2 $ vertex pairs
$\{u,v\}$, and the co-degree of each $\{u,v\}$ is the
sum of $(N-2, \floor {\frac k 2}, p^2)$-variables. \qed
\fi
%\ifPersonal Proof

%========================================================
%                   Jumbledness
%========================================================
The following definition is a modified version of the one in
\cite{t,cgw}, see also \cite{ac} and \cite{as}, Chapter 9.

\bdf {Jumbledness} \label{def_jumbl}
For vertex sets $U,V$, let $e(U,V)$ denote the number of edges that
pass from $U$ to $V$ (internal edges of $U \bigcap V$
are counted twice). A graph is $(p,\alpha)$-jumbled if
$e(U,V)=p |U||V| \pm \alpha \sqrt{|U||V|}$ holds for all
$U,V$. \edf

\btn{Achieving optimal jumbledness}
\label{kwise_gives_jumbl} There exist absolute
constants $c_1,c_2,c_3$ \st all \kwig \engnk are \as
$(p,\alpha)$-jumbled whenever either:
\ben
\item
$k\geq 4$, $p\geq \Omega(\frac {1}{N})$~~and~~$\alpha
\gg \sqrt{p}N^{3/4}$;~~or
\item
$k\geq \log(N)$, $\frac {c_1 \log(N)}{N} \leq p \leq
1- \frac {c_2 \log^4(N)}{N}$~~and~~$\alpha \geq c_3
\sqrt{pN}.$ \een \etn
\noindent{\bf Proof.} The proof is based on spectral techniques
and combines some refined versions of ideas  from \cite{ac}, \cite{fk2}
and \cite{vu}, using the fact that  traces of the $k$-th power of the
adjacency matrix of a graph are identical in the $k$-wise independent case
and in the totally random one. The details are somewhat lengthy and
are thus deferred to Appendix \ref{proof_kwise_gives_jumbl}.

\mysubsection{Connectivity, Hamiltonicity and Perfect
Matchings}\label{section_conn_Ham}

%=======================================================================
%                   Connectivity
%=======================================================================

\btn{Achieving connectivity}
\label{kwise_gives_conn}
There exists a constant $c$ \st the following holds. All \kwig \engnk are \as
connected whenever either:
\begin{itemize}
\item
$k\geq 4$~~and~~$p \gg \frac {1}{\sqrt{N}}$;~~or
\item
$k\geq 4\log(N)$~~and~~$p\geq \frac {c\ln(N)}{N}$.
\end{itemize}
\etn
\noindent{\bf Proof.} Let $U$ be a vertex-set that
induces a connected component. Connectivity follows from
having $|U|>0.5N$ for all such $U$. The following holds
\as for \gnkns. By Lemma \ref{kwise_gives_deg}, all
vertices have degree $\geq 0.9pN$, so $e(U) \geq
0.9pN|U|$. By Theorem \ref{kwise_gives_jumbl}, all sets
$U$ satisfy $e(U)\leq p|U|^2 + \alpha|U|$ with $\alpha=
O(\sqrt{pN})= o(pN)$. Re-arranging gives $(0.9-o(1))N
\leq |U|$. \qed

%=======================================================================
%                   Hamiltonicity
%=======================================================================

\btn{Achieving Hamiltonicity}
\label{kwise_gives_Ham}
All \kwig \engnk are \as Hamiltonian (and for even $N$ contain a perfect matching) whenever either:
\begin{itemize}
\item
$k\geq 4$~~and~~$p \geq \frac {\log^2(N)} {\sqrt
N}$;~~or
\item
$k\geq 4\log(N)$~~and~~$p\geq \frac {\log^2(N)}{N}$.
\end{itemize}
\etn
\noindent{\bf Proof.} Let $\bar \Gamma(V)$ denote the
set of vertices $v \notin V$ that are adjacent to some
vertex in the vertex-set $V$. By Theorem 1.1 in
Hefetz, Krivelevich \andCoAuthers Szabo's \cite{hks},
Hamiltonicity follows from the existence of constants $b,c$
such that \as
(i) $|\bar \Gamma(V)|\geq 12 |V|$ holds for all sets
$V$ of size $\leq b N$, and (ii) $e(U,V)\geq 1$ holds
for all disjoint sets $U,V$ of size $\frac {c
N}{\log(N)}$. We remark that (unlike other asymptotic
arguments in this paper), the sufficiency of (i) and
(ii) might hold only for very large $N$.
For (i), let $b= \frac 1 {170}$ and consider an
arbitrary set $V$. By Theorem \ref{kwise_gives_jumbl},
\as all vertex-sets $T$ have $e(T) \leq p|T|^2 +
o(pN)|T|$. By Lemma \ref{kwise_gives_deg} \as all the
degrees are $\opmo pN$, so exactly $\opmo pN|V|$ edges
touch $V$ (where internal edges are counted twice).
Let $T=V \bigcup \bar \Gamma(V)$, and assume that $|\bar
\Gamma(V)| < 12|V|$. We get $\omo pN|V| \leq e(T) \leq p
(13|V|)^2 + o(pN)|V|$. Re-arranging gives $|V| > \frac N
{170}$. Condition (i) follows.
For (ii), by Theorem \ref{kwise_gives_jumbl}, \as all
(equal-sized and disjoint) vertex-sets $U,V$ have
$e(U,V) \geq p|U||V| - O(\sqrt{pN})|U|$. If there is
no edge between $U$ and $V$, then $e(U,V)=0$.
Re-arranging gives $|U|\leq O(\sqrt{N/p}) \leq O(\frac
N{\log(N)})$. Condition (ii) follows. \qed

%=======================================================================
%            Connectivity - Negative Result
%=======================================================================
%
\btn {Failing to preserve connectivity}\label{kwise_gives_Ham} There exist
pair-wise independent graphs \engnk\\
where $p=1/2$, that are (i) \as disconnected
(and contain no Hamiltonian cycles), and that (ii)
contain no perfect matchings \wip  $1$. \etn
\noindent{\bf Proof.} Consider the graphs defined by
partitioning all vertices into 2 disjoint sets
$V_0,V_1$ where each $V_j$ induces a clique, no edges
connect $V_0$ to $V_1$, and $V_1$ is chosen randomly and uniformly among
all subsets of odd cardinality  of the vertex set. Note that for every
set of $4$ vertices, there are $16$ ways to split its vertices among $V_0$
and $V_1$, and it is not difficult
to check that if $N \geq 5$, then each of these $16$ possibilities
is equally likely.
Therefore, any edge appears w.p.~$\frac 1 2$, and any pair
of edges (whether they share a common vertex or not)
appears w.p.~$\frac 1 4$. Still the graph is connected
iff all the vertices belong to the same $V_j$ which
happens only w.p.~$2^{-N+1}$ (and only if $N$ is odd).  Since
$|V_1|$ is odd, the graph contains no perfect matching.
\qed
%
%------------------------------------------------------
%                    3-wise comment
%------------------------------------------------------
%
\ifPersonal %
Note that these graphs are not 3-wise independent, as
the existence of the edges $\{u,v\}$ and $\{v,w\}$
implies the existence of the edge $\{u,w\}$. Also note
that this construction cannot be generalized for $p \neq
\frac 1 2$ (by assigning vertices w.p.~$q$ to $V_0$ and
w.p.~$1-q$ to $V_1$). Indeed, consider the requirements
(i) $p=q^2+(1-q)^2$ and (ii) $p^2=q^3+(1-q)^3$. If we
plug (i) into (ii) we get an equation of degree 4 that
by Vieta's formula has only 3 solutions $p=0,\frac 1
2,1$ (see
\url{http://planetmath.org/encyclopedia/VietasFormula.html}).

We also remark that when $p$ is slightly increased to
$1/2 + N^{-\Theta(1)}$, then 4-\wice {\em does} suffice
for achieving connectivity since then all vertices \as
have degree $> N/2$. Then, by Dirac's Theorem, the
graphs are in fact Hamiltonian. \else
%
%------------------------------------------------------

Note that when $p$ is slightly increased to $1/2 +
N^{-\Theta(1)}$, then 4-\wice suffices for achieving
Hamiltonicity (via Dirac's Theorem), because then \as
all vertices have degree $> N/2$.
\fi%

%===================================================
%                  Kappa via Jumbledness
%===================================================

\mysubsection{High-connectivity}
\label{section_high_conn}
\btn{Achieving optimal connectivity}
\label{kwise_gives_kappa_jumbl}
There exists an absolute constant $c$, \st for all
\kwig \engnk the connectivity number is \as $\opmo pN$
when either
\begin{itemize}
\item
$k\geq 4$~~and~~$p \gg N^{-\frac 1 {3}}$;~~or
\item
$k\geq \log(N)$~~and~~$p \geq c \frac{\log(N)}{N}$.
\end{itemize}
\etn
\noindent{\bf Proof.}
The connectivity is certainly not larger than $(1+o(1))pN$, as it
is upper-bounded by the minimum degree.
By Theorem 2.5 in Thomason's \cite{t} $\kappa \geq d-
\alpha/p$ holds for any $(p,\alpha)$-jumbled graph with
minimal degree $\geq d$. Thus, achieving $\kappa \gtrsim
pN$, reduces to obtaining (i) $d = (1\pm o(1))pN$, and
(ii) $\alpha \ll pd$.
Condition (i) \as holds by Lemma \ref{kwise_gives_deg}.
By Theorem \ref{kwise_gives_jumbl}, we \as achieve
$(c_3\sqrt{pN})$-jumbledness for some constant $c_3$,
so condition (ii) becomes $p^2N \gg \sqrt{pN}$. This
proves the first part of the theorem. To prove the
second we note, first, that we may assume that $p \ll
1$ (since otherwise $4$-wise independence suffices).
Let $S$ be a smallest separating set of vertices,
assume that $|S|$ is smaller than $(1-o(1))pN$, let
$U$ be the smallest connected component of $G-S$ and
let $W$ be the set of all vertices but those in $U
\cup S$. Clearly $|W| \geq (\frac{1}{2}-o(1))N$. Note
that $e(U,W)=0$, but by jumbledness $e(U,W) \geq
p|U||W| -c_3 \sqrt {pN |U||W|}$. This implies, using
the fact that $|W|>N/3$, that $|U| \leq
\frac{3c_3^2}{p}$. Using jumbledness again, $e(U,S)
\leq p|U||S|+c_3 \sqrt {pN |U||S|}$ but as all degrees
are at least $(1-o(1))pN$, $e(U,S) \geq
(1-o(1))pN|S|-e(U) \geq (1-o(1))pN|U|-p|U|^2- c_3
\sqrt {pN} |U| \geq |U| (1-o(1))pN$, where here we
used the fact that $|U| \leq O(1/p)$ and that
$\sqrt{pN} =o(pN)$. This implies that either $p|U||S|
\geq \frac{1}{2}|U|pN$, implying that $|S| \geq N/2
\gg pN$, as needed, or $c_3 \sqrt {pN|U||S|}  \geq
\frac{1}{3}|U|pN$, implying that $|S| \geq
\frac{1}{9c_3^2} |U| pN$ which is bigger than $pN$
provided $|U| \geq 9c_3^2$. However, if $|U|$ is
smaller, then surely $|S| \geq (1-o(1))pN$, since all
degrees are at least $(1-o(1))pN$ and every vertex in
$U$ has all its neighbors in $U \cup S$.
\qed

%===================================================================
%        {Preserving the Threshold for Subgraphs}
%===================================================================
%
\mysubsection{Thresholds for the Appearance of
Subgraphs}\label{section_sub_graphs}

For a fixed non-empty graph $H$, let $\rho(H)$ and
$p^*_H$ be as in Section \ref{intro_H_copies}.

\bob {Preserving the threshold for appearance of
sub-graphs} \label{kwise_gives_subgraph_threshold}
There exists a function $D(v)=\opmo \frac {v^4}{16}$
\st for any graph $H$ with at most $v$ vertices, and
for all \kwig \engnk $ $ with $k \geq D(v)$ the
following holds. Let $A$ denote the event that $H$
appears in \gnk (not necessarily as an induced
sub-graph). Then
\begin{itemize}
\item
If $p(N) \ll p^*_H(N)$ then $(\neg A)$ \as holds.
\item
If $p(N) \gg p^*_H(N)$ then $A$ \as holds.
\end{itemize} \eob
\noindent{\bf Proof.}
The proof (given in Appendix
\ref{kwise_gives_subgraph_thresholdPrf}) applies
Rucinski \andCoAuthers Vince's \cite{rv} to derive a
lower-bound on the minimal $k$ sufficient for the
original \gnp argument to hold. \qed

%===================================================================
%         {Defying the Threshold for Subgraphs}
%===================================================================
%
\btn {Defying the threshold for appearance of sub-graphs}\label{sub_graphs_defy}
For any (fixed) graph $H$ that satisfies%
\footnote{This condition rules out only graphs $H$ that
are a collection of disjoint edges. For such graphs
$\rho(H)=2$, so clearly no $H$-copies can be produced
(even if $k=1$) when $p(N) \ll p^*_H(N)= N^{-2}$.}
$\rho(H)<2$, there exists \kwig \engnk where
$k=\ceil{\frac 2 {\rho(H)} -1}$ and $p(N) \ll p^*_H(N)$
\st $H$ \as appears in \gnk as an induced sub-graph.
\etn

\ifPersonal %-----------------------------------------
\begin{remark}
The term $\ceil{\frac 2 {\rho(H)} -1}$ is merely the
minimal $k$ \st $k \gneq \frac 2 {\rho}$ (verified by
separately checking for either integer or non integer
$\frac 2 {\rho}$).\end{remark}

\begin{remark} \label{remark_rho_cond_reasonable}
Note that whenever $\rho \lneq 2$, we have
$\ceil{\frac 2 {\rho(H)} -1} \geq 1$, so the condition
on $k$ is reasonable.\end{remark}
\fi % \ifPersonal remark
%
%-----------------------------------------------------

\noindent{\bf Proof.} Theorem \ref{sub_graphs_defy}
relies on Lemma
\ref{kwise_construction_unexpected_sub_graphs}. This
lemma considers the appearance of the sub-graph $H_N$ in
\gnk where $\ensm{H_N}$ is any sequence of graphs
(possibly) with unbounded order.

%=======================================================================
%               Unexpected  Sub-Graphs
%=======================================================================
%
\ble {\kwig with unexpected appearance of sub-graphs}
\label{kwise_construction_unexpected_sub_graphs}
Let $\ensm{H_N}$ be a sequence of graphs where $H_N$ has
exactly $S(N) < \sqrt N$ vertices, $e_1(N)$ edges and
$e_0(N)$ none-edges. Assume that for each $N$ there
exists $(\binom {S(N)} 2,k(N),p(N))$-variables \st with
probability $\Delta(N) \gg (S(N)/N)^2$ it holds that
the first $e_0(N)$  variables attain value $0$
and the next $e_1(N)$ variables attain value
$1$. Then there exist \kwig \engnk that \as contain
$H_N$-copies as induced sub-graphs. \ele
%----------------------------------------------
%
\noindent{\bf Proof (Lemma
\ref{kwise_construction_unexpected_sub_graphs}).} Fix
$N$, so $H=H_N, S=S(N), e_i=e_i(N), k=k(N), p=p(N),
\Delta=\Delta(N).$ We construct graphs \gnk that \as
contain $H$ copies. Given the $N$ vertices,
let $\{V_j\}_{j=1}^M$ be any maximal collection of {\em
edge-disjoint} vertex-sets, each of size $|V_j|=S$.
For each $j$, decide the internal edges of $V_j$ by some $(\binom S 2,k,p)$-variables \st $H$ is induced by $V_j$ with probability $\Delta$. This can be done by appropriately defining which specific edge in $V_j$ is decided by which specific variable.
Critically, the constructions for distinct sets $V_j$
are totally independent. The $R=\binom N 2 - M\binom S
2$ remaining edges
can be decided by any $(R,k,p)$-variables. The
resulting graph is clearly \kwins.

The main point is that (i) the events of avoiding
$H$-copies on the various sets $V_j$ are totally
independent (by the edge-disjointness of the $V_j$-s),
and that (ii) in our \kw case $\Delta$ is rather large
(compared with the totally independent case). Thus,
avoiding $H$-copies on {\em any} of the $V_j$-s is
unlikely.
Indeed, let $B$ denote the event that no $H$-copies appear in the resulting graph, while $B'$ only denotes the event that none of the $V_j$-s induces $H$.
By Wilson's \cite{wil} and Kuzjurin's \cite{kuz} we have $M= \Theta(N^2/S^2)$\ifPersonal~(which is obviously optimal), \else, \fi so

$$\pr{B} \leq \pr{B'}
= (1-\Delta)^{M} \leq
e^{-\Theta\lb\frac {\Delta N^2}{S^2}\rb},$$
which vanishes by our requirement that
$\Delta \gg (S/N)^2.$ %
\qed %
~(Lemma \ref{kwise_construction_unexpected_sub_graphs})

\noindent{\bf Completing the proof of Theorem
\ref{sub_graphs_defy}}.
For $v=v(H), \rho=\rho(H), p^*=p^*_H$,
and some $1 \ll f(N)\leq N^{o(1)}$, define $p$ \st $p^{-1}$ is the minimal power of 2 that is larger than, $\frac {f(N)}{p^*}$. As desired $p \ll p^*$.
Let $e_1$ and $e_0$ respectively denote the number of edges and non-edges in $H$.
With $M=\binom v 2$ and $F=1/p$, we apply Lemma \ref{modify_JoffeCG_lem} to produce $(M,k,p)$-variables
\st \wip $\geq F^{-k}$ the first $e_0$ variables have value 0, and the remaining $e_1$ variables have value 1.
By Lemma \ref{kwise_construction_unexpected_sub_graphs}, the latter immediately implies the existence of \kwig that \as contain $H$-copies as long as $F^{k} \ll (N/v)^2$. As $F=1/p=N^{\rho+o(1)}$, this $\ll$ requirement translates to $k\rho \lneq 2$. \qed
~(Theorem \ref{sub_graphs_defy})

%=======================================================================
%               Chromatic Number
%=======================================================================

\mysubsection{The Chromatic Number}\label{section_color}

%======================================================
\bob {Preserving the chromatic-number lower bound}
%======================================================
%
\label{chi_lo_bound} For any $c>0$ there exists some
$d>0$, \st all \kwig \engnk with $(\log(N))^{-c} \leq
p \leq 1 -N^{-o(1)}$ and $k \geq d (\log(N))^{c+1}$
\as have chromatic number $\chi \geq \frac
{N\log(1/1-p)} {2\log(pN)}$. \eob

\ifPersonal
\noindent{\bf Proof.} Let $I(G)$ denote the independence
number of (a single) $N$-vertex graph $G$, so
$\chi(G) \geq \frac N {I(G)}$.
Whenever $k \geq \Theta \lb (\frac 1 p \log(pN))^2
\rb$ our \kwice upper-bound on I (observation
\ref{kwise_gives_exact_indp_num}) shows that \as
$I(\gnkm) \leq \frac {2\log(pN)}{\log(1/1-p)}= O\lb
\frac 1 p \log(pN) \rb$ (the $=O(\cdot)$ applies
$\log(1/1-p) \stackrel {p \rightarrow
0}{\longrightarrow} p$),
so $\chi(\gnkm) \geq \Omega \lb \frac
{pN}{\log(pN)}\rb$. \qed \else \noindent{\bf Proof.}
Let $I(G)$ denote the independence number of (a
single) $N$-vertex graph $G$. Clearly, $\chi(G) \geq
\frac N {I(G)}$, so observation \ref{chi_lo_bound}
follows from our \kwice upper-bound on I (observation
\ref{kwise_gives_exact_indp_num}). \qed \fi

%======================================================
\btn {Preserving the chromatic-number upper bound}
%======================================================
%
\label{chi_up_bound} There exists an absolute constant
$c$ \st the following holds. All \kwig \engnk with $p
\leq 1/2$ \as have chromatic number $\chi \leq \frac
{cN\log(1/1-p)} {\log(pN)}$, whenever either: \ben
\item
$k\geq 12$~and~$p \geq N^{-\frac 1 {75}}$;~~or
\item
$k\geq \log(N)$~and~$p \geq c \frac {\log(N)}{N}$.\een
\etn

\noindent {\bf Remark.} No special effort was made to
optimize the constants $\frac 1 2$ and $\frac 1 {75}$.

\noindent{\bf Proof (sketch).} Since $p$ is bounded from above
and $\log(1/1-p) \stackrel {p \rightarrow 0} {
\longrightarrow} p/\ln(2)$, it suffices to show that
\as $\chi \leq  O ( \frac {pN} {\log(pN)} )$.
Item 1 is based on Alon, Krivelevich \andCoAuthers
Sudakov's \cite{aks}.
Specifically, choose $\delta=1/25$, \st by item 1 in
Lemma \ref{kwise_gives_deg} (with $\eps=(\log(N))
p^{-1/2} N^{-3/8}$) and by item 1 in Lemma
\ref{kwise_gives_codeg} (with $\gamma= (\log(N))
p^{-1} N^{-1/6}$), \as all the degrees are lower
bounded by $pN(1- p^{-1/2} N^{-3/8+o(1)}) \geq pN-
N^{1-4\delta},$ and all co-degrees are upper bounded
by $p^2N (1+ p^{-1} N^{-1/6+o(1)}) \leq p^2N-
N^{1-4\delta}.$
By Theorem 1.2 in \cite{aks}, these conditions (with
$\delta <1/4$ and $p\geq N^{-\frac {\delta}{3}}$)
imply that $\chi \leq \frac {4pN}{\delta \ln N} \leq
O(\frac{pN}{\log(pN)}).$

Item 2 follows from jumbledness and the main result of
\cite{aks1} (which is based on \cite{joh}), by which
any graph with maximum degree $d$ in which every
neighborhood of a vertex contains at most
$d^{2-\beta}$ edges (for some constant $\beta$) has
chromatic number $\chi \leq O(\frac d {\log d})$. \qed

%==================================
\paragraph{\Large Acknowledgements}
%==================================
The second author wishes to thank Oded Goldreich for
his encouragement, and Ori Gurel-Gurevich, Chandan
Kumar Dubey, Ronen Gradwohl, Moni Naor, Eran Ofek, Ron
Peled, and Ariel Yadin for useful discussions.

%=======================================================================
%   CITATIONS
%=======================================================================

%=======================================================================
%                     Appendix
%=======================================================================
%
%=======================================================================
%         Modified Construction (M,k,p) Variables
%=======================================================================
\mysection{Appendix - Detailed Proofs}
\mysubsection{Modified Construction of ${\bf k}$-Wise
Independent Variables - Proving
Lemma~\ref{modify_JoffeCG_lem}}
\label{modify_JoffeCG_sec}
Recall that given any prime power $F$, the original
\cite{jof,cg,abi} construction considers the field $\mathbb F$ with elements $\{0,...,F-1\}$, and for each element $j \in \mathbb F$, a random variable $Z_j$ is defined, \st the $Z_j$s are $k$-wise independent, and each $Z_j$ is uniformly distributed in $\{0,...,F-1\}$.
We derive from those $Z_j$-s some $(M,k,p)$ binary variables $X_j$, by setting %
(i) $X_j=1$ iff $\frac {Z_j+1}{F} \geq 1-p$ for $j=1,...,e_0$, and
(ii) $X_j=1$ iff $\frac {Z_j+1}{F} \leq p$ for $j=e_0+1,...,e_0+e_1$.
Evidently, the $X_j$-s are \kwi with $\Pr(X_j=1)=p$.
\ifPersonal Thus, the desired event is $A=\lb
\bigwedge_{j=1}^{e_0} (X_j=0) \rb
\bigwedge \lb
\bigwedge_{j=e_0+1}^{e_0+e_1} (X_j=1) \rb$ %
\fi %\ifPersonal
Recall that $Z_j \eqdef Q(j)$ with $Q$ being a uniformly
random degree $k$ polynomial over $F$, and let $B$
denote the event that the 0-polynomial was chosen. Since
$B$ implies $A$, we get $\Pr[A] \geq \Pr[B]= F^{-k}$. \qed

\ifPersonal \begin{remark} If we considered exactly
$k$ edges we could improve $\Pr[A] \geq p^{k}$.
Indeed, fix any $k$ field elements $x_1,...,x_k$, and
consider the interpolation-bijection between the set
of degree $k$ polynomials $Q$, and the set of vectors
$Q(\vec x)=(Q(x_1),...,Q(x_k))$. Clearly, $A$ holds
for precisely a $p^k$ fraction of the $Q(\vec x)$, and
hence for the same fraction of the polynomials $Q$,
giving $\Pr[A] \equiv p^{k}$.
The point is that we consider more than $k$ edges.
Thus, when the first $k$ variables receive the desired
values (an event that {\em does} hold \wip $\equiv
p^k$), the interpolation may always provide the
undesired valued for the remaining variables. This
unfortunate outcome is ruled out only in the specific
case of the 0-polynomial.\end{remark}
\fi %complete remark

%=======================================================================
%            k-Wise Chebyshev-Bound Proof
%=======================================================================
\mysubsection{${\bf k}$-Wise Independence Tail Bound -
Proving Lemma
\ref{kwise_Chebyshev_bound}}\label{sec_inequal}
Let $\bar X_i  \eqdef  X_i - \mu$ and $\bar X \eqdef
\sum_{i=1}^{{M}} \bar X_i$, so $ X-\mathbb E (X)=\bar
X.$ Thus,

\begin{eqnarray}
\Pr[|X-\mathbb E (X)| \geq  \delta \mathbb E (X)] &=& \Pr[|\bar X| \geq \delta \mathbb E (X)] \nonumber\\
= \Pr[\bar X^{k} \geq (\delta \mathbb E (X))^{k}]
&\leq& \frac {\mathbb E (\bar X^{k})}{(\delta \mathbb E (X))^{k}}, \nonumber
\end{eqnarray}
the last equality holds for any even positive $k$, while the $\leq$ employs Markov's inequality.

We bound $E  \eqdef  \mathbb E (\bar X^{k})$ using the
expansion
$$\bar X^{k} =
\sum_{\vec{d} \in D} \Pi_{i=1}^{{M}} \bar X_i^{d_i},$$
where $D  \eqdef  \{\vec{d}=(d_1,...,d_{{M}}) | \sum
_{i=1}^{{M}} d_i=k, d_i\geq 0\}$. The $k$-wise
independence of the variables $X_i$, guarantees the
$k$-wise independence of the $\bar X_i$s, so
$$E = \sum_{\vec{d} \in D} \mathbb E (\Pi_{i=1}^{{M}} \bar X_i^{d_i})=
\sum_{\vec{d} \in D} \Pi_{i=1}^{{M}} \mathbb E (\bar
X_i^{d_i}).$$
Next, since $\mathbb E(\bar X_i)=0$ we can ignore all
terms where $d_i=1$ for some $i$. Namely, we consider only terms $\Pi = \Pi_{\ell=1}^{j} \mathbb E (\bar
X_{i_{\ell}}^{d_{i_{\ell}}})$ where for some $j \leq
\fkt$, it holds that precisely $j$ variables appear and for each variable $\bar X_{i_{\ell}}$ in $\Pi$ we have $d_{i_{\ell}} \geq 2$.
Hence, \begin{equation}E \leq \sum_{j=1}^{\fktb} \Psi_j \label{E_def},\end{equation}
whenever $\Psi_j$ bounds the contribution of all terms
$\Pi$ with precisely $j$ variables.

Strengthening standard versions of the inequality begins by taking
\begin{equation}\label{psi_def}
\Psi_j \eqdef { M \choose j} j^{k}
\lsb\mu(1-\mu)\rsb^j.\end{equation}
Indeed, ${ M \choose j} j^{k}$ clearly bounds the
number of terms $\Pi$ with precisely $j$ variables,
while $\lsb \mu(1-\mu)\rsb^j$ bounds the expectation of each term $\Pi$ that has precisely $j$ variables because
\begin{eqnarray}
\mathbb E [(\bar X_i)^d] &=& \mu(1-\mu)^d +(1-\mu)(-\mu)^d \nonumber\\
&\leq& \mu(1-\mu)^d +(1-\mu)(+\mu)^d \nonumber\\
&\leq& \mu(1-\mu)[(1-\mu)^{1} +\mu^{1}]
=\mu(1-\mu) \label{exp_X_i_to_d_bound}
\end{eqnarray}
(the final $\leq$ applies the facts $0 \leq \mu,1-\mu \leq 1$ and $d \geq2$).
Thus, multiplying over the $j$ terms gives
$$\Pi_{{\ell}=1}^{j} \mathbb E (\bar
X_{i_{\ell}}^{d_{i_{\ell}}}) \leq \lsb \mu(1-\mu)\rsb^j.$$

Observe that $\Psi_j$ is maximized when $j=\fkt$.
Indeed,
\begin{eqnarray}
\frac {\Psi_{j+1}}{\Psi_j} &=& \frac {M-j}{j+1}\left(\frac{j+1}{j}\right)^k \mu(1-\mu) \nonumber\\
&>& \frac {M-k}k \mu(1-\mu) \geq 1 \nonumber
\end{eqnarray}
(the concluding $\geq 1$ holds by the lemma's assumption).

Thus, the maximal $\Psi_j$ is
\begin{eqnarray}
\Psi_{\hspace{-0.4ex}\fktb} &=& {M \choose {\fkt}} \fktbb^k [\mu(1-\mu)]^{\fktb} \nonumber\\
&\leq& \frac {M^{\fktb}} {(\fkt)!} {\fktbb}^k [\mu(1-\mu)]^{\fktb} \nonumber\\
&\leq& \frac {(e M)^{\fktb}} {\sqrt{2\pi \fktb}(\fkt)^{\fktb}} {\fktbb}^k [\mu(1-\mu)]^{\fktb} \nonumber\\
&=& \frac {\lsb\frac e 2 M k \mu(1-\mu)\rsb^{\fktb}} {\sqrt{\pi k}} \nonumber
\end{eqnarray}
(Stirling's approximation for $(\fkt)!$ implies the last
$\leq$).

To summarize, all the above gives
\begin{eqnarray}
\Pr \lsb|X-\mathbb E (X)| \geq \delta \mathbb E (X)\rsb & \leq & \frac {k \Psi_{(k/2)}}{\left (\delta \mathbb E(X) \right)^{k}}\nonumber\\
& \leq & \frac {\frac k {\sqrt{\pi k}} (\frac e 2 M k \mu(1-\mu))^{\fktb}}{\left (\delta  \mu M \right)^{k}} \nonumber\\
& \leq & \sqrt{\frac k {\pi}} \left [ \frac {\frac e 2 k(1-\mu)}{\delta^2 \mu M} \right]^{\fktb}.\nonumber
\end{eqnarray}
The Lemma follows as it can be directly shown that for
all $k$
$$\sqrt{ \frac k {\pi}} \lb\frac e 2\rb^{\fktb} \leq 2^{\fktb}. \qedm$$

%=======================================================
%   Tightness discussion - Removed to OMITTED_K_WISE.tex
%=======================================================

%=======================================================================
%        Match Lemma - Removed to OMITTED_K_WISE.tex
%=======================================================================

%=======================================================================
%             Sub-Graphs - Preserve Proof
%=======================================================================

\mysubsection{Appearance of Subgraphs - Proving
observation \ref{kwise_gives_subgraph_threshold}}
\label{kwise_gives_subgraph_thresholdPrf}

We first consider only balanced graphs $H$, namely graphs where $\rho(H) \leq \rho(H')$ for any subgraph $H' \subseteq H$. The original
\gnpns-threshold proof \cite{erGiant} takes a fixed
graph $F$ as a parameter, and considers for each set
$T$ of $v(F)$ distinct vertices the random variable
$Y^{F}_T$ which indicates whether $T$ spans $F$ in the
resulting graph. Thus $Y^{F}\eqdef \sum_{T}Y^{F}_T$
counts the number of sets that span $F$.

First, the authors of \cite{erGiant} consider a specific subgraph $H'
\subseteq H$ \st $\rho(H)=\frac {v(H')} {e(H')}$ and
show that $p \ll p^*_H$ implies that $\mathbb E(Y^{H'})
\ll 1$. In this case, $H'$ rarely appears in \gnp graphs
and so does $H$. On the other hand, whenever $p \gg
p^*_H$, they show that $\mathbb E(Y^{H}) \gg 1$ and by
Chebyshev's inequality it is deduced (only here the fact
that $H$ is balanced is
used), that %
\ifPersonal \as the number, $Y^{H}$, of vertex sets that span $H$-copies is close to its expectation ($Y= \opmo \mathbb E (Y)$) and in particular \fi
some $H$-copies appear. Thus, the entire argument
applies only probabilities regarding either a
single variable $Y^{F}_T$, or a pair
$Y^{F}_T,Y^{F}_{T'}$ of variables, and relies only
upon the independence of sets of $m \leq 2 {{v(H)} \choose 2}$ edges.

For non-balanced graphs the $p \ll p^*_H$ part holds
as for balanced ones. For $p \gg p^*_H$, we rely on
the fact that for any graph $H$, there exists an
extension graph $H \subseteq H''$ s.t. $H''$ is
balanced and $\rho(H'')=\rho(H)$ (Rucinski
\andCoAuthers Vince
\cite{rv}). %
Since $p \gg N^{-\rho(H)}$ means that $p \gg N^{-\rho(H'')}$, and since $H''$ is balanced, then
\gnp graphs \as contain copies of $H''$, and copies of
$H$ appear as well. This time the Chebyshev argument
assumes only the independence of sets of
$m \leq 2 {{v(H'')} \choose 2}$ edges. Since by \cite{rv} there
exists  $H''$ as above with $v(H'') \leq \opo \frac {\lsb v(H)\rsb ^2} 4$,
then $m=\opmo \frac {\lsb v(H)\rsb^4}{16}$ suffices. \qed

%=======================================================================
%            Appendix - Independence Number
%=======================================================================

\mysection{Appendix - The Independence Number of ${\bf
k}$-Wise Independent graphs}\label{appendix_indp_num}

%=======================================================================
%            Precise Independence Number
%=======================================================================

The following positive result follows the argument used
to establish observation
\ref{kwise_gives_subgraph_threshold}.

\bob {Preserving random graphs' precise
independence-number} \label{kwise_gives_exact_indp_num}
Consider arbitrary \kwig \engnk where $N^{-o(1)} \leq
p(N) \leq 1-N^{-o(1)}$, and let $I(N)=I(\gnkm)$ denote
the independence number of \gnkns.
Then there exists a function $S^*(N,p)=\omo \frac
{2\log(pN)} {\log(1/(1-p))}$, \st if $k(N)\geq
S^*(N,p)+2$, then \as $I(N)\leq S^*(N,p)+1$, and if
$k(N)\geq \binom {S^*(N,p)} 2$, then \as $I(N)\geq
S^*(N,p)$. \eob
%----------------------------------------
\noindent{\bf Proof.}
The classical proof (\cite{be}, \cite{mat}) of this claim for \gnp
graphs considers for each set $T$ of $S$ distinct vertices ($S$
being a parameter) the random variable $Y^S_T$ which indicates
whether $T$ spans an independent set in the resulting graph.
Thus $Y^S=\sum_T Y^S_T$ counts the total number of independent sets of size $S$.
It is shown that for $S=S^*+2$ then $\mathbb E(Y^S)\ll 1$ so \as the independence number $\leq S^*+1$. On the other hand, for $S=S^*$ then $\mathbb E(Y^S)\gg (1)$, and by Chebyshev's inequality it is deduced that \as \ifPersonal the number of independent sets of size $S^*$, $Y^S$, is close to its expectation and in particular \fi some independent sets of size $S^*$ appear.
This entire argument considers only probabilities regarding either a single variable $Y^S_T$ (for the lower- and upper-bound on $I$), or a pair $Y^S_T,Y^S_{T'}$ of variables (for the lower-bound). Therefore, the upper-bound holds for all \kwig with $k \geq S^*+2$, and the lower-bound holds whenever $k \geq 2 {{S^*} \choose 2}$.
\ifPersonal%
\footnote{Indeed, let $S=\frac {2\log_2N}
{\log_2(1/(1-p))}+1$. The expected number of
independent sets of size $S$ is $\binom N S
(1-p)^{\binom S 2} < \lsb \frac {eN}S (\frac 1
{(1-p)})^{(-\log(N)) / (\log(1/(1-p)))}\rsb^S= [\frac
{eN}S 2^{-\log(N)}]^S= (\frac e S)^S$, which clearly
vanishes whenever $S\gg 1$. Finally, $S\gg 1$ holds
iff $p\leq 1-N^{-o(1)}$ and the claim follows.}
\else~\fi \qed
%=======================================================================
%              Huge Independence Number
%=======================================================================

We next provide our negative results. Since the complexity of known constructions of \kwv critically depend on the length, $\ell(p)$, of the binary representation of $p=0.b_1...b_{\ell}$, it is reasonable to focus on densities with bounded length.
The argument used here was already applied in the context of Theorem \ref{sub_graphs_defy}.

\btn {K-wise independent graphs with huge independent sets} \label{kwise_with_huge_indp_set}
Let $S,k:\N \into \N^+$ and $p:\N \into (0,1)$ satisfy
$S(N)\ll N^{\lb \frac 1 {k(N)+1} \rb }$ and
$\ell(p(N)) \leq 2 \log \lb S(N)\rb$. Then there exist
\kwig $ $ \engnk that \as contain independent sets of
size $S(N)$. \etn

%=======================================================================
%      Huge Independence Number - Theorem Proof
%=======================================================================

\noindent{\bf Proof.} By Lemma
\ref{modify_JoffeCG_lem}, since $\ell(p) \leq 2 \log
S$, we get $(\binom S 2,k,p)$-variables \st the
probability that all \ifPersonal the binary \fi
variables receive value 0 is $\Delta \geq S^{-2k}$.
 From this, Lemma \ref{kwise_construction_unexpected_sub_graphs}
gives \kwig that \as contain independent-sets of size $S$, whenever $S^2 \ll \Delta N^2$ \ifPersonal
which is guaranteed by $S^{2k+2} \ll N^2$. \else. \fi
\qed

%=======================================================================
%        Corollary for Huge Independence Number
%=======================================================================
%
\bcorl \label{kwise_with_huge_indp_set_corol}
Let $(S,k,p)$ be as in Theorem \ref{kwise_with_huge_indp_set}, with $\Omega(1) \leq p(N)\leq 1-N^{-o(1)}$, and with $S^*$ as in observation \ref{kwise_gives_exact_indp_num}%
\ifPersonal denoting the \gnp \as upper-bound on the
independence number. \else. \fi Fix $c>1$.
Then there exist \kwig \engnk where $k(N)\geq \omo
\frac {\log(N)}{c\log \log(N)}$ that \as contain
independent sets of size $S \gg \lb S^*(N)\rb^{c}$.
\ecorl

%=======================================================================
%      Huge Independence Number - Prove Corollary
%=======================================================================

\noindent{\bf Proof.}
It suffices to provide an integer $S$ \st: (i) $S \gg \lb S^* \rb^c$ (the desired outcome) and (ii) $S \ll N^{\lb \frac 1 {k+1} \rb }$ (the sufficient condition for applying Theorem \ref{kwise_with_huge_indp_set}). Such $S$ clearly exists as long as
\begin{equation}\label{eq__1}
\lb S^* \rb^c \ll N^{\lb \frac 1 {k+1} \rb }.
\end{equation}
Define $r$ by $S^*=N^{r}.$ Since $N^{\frac
{1}{\log(N)}}=2$, then any choice of $f(N)\gg 1$ gives
$N^{\frac {f(N)}{\log(N)}} \gg 1.$ Thus, equation
(\ref{eq__1}) translates to having (for some $f(N)\gg
1$)
\begin{equation}\label{eq__2}
cr \leq  \frac 1 {k+1} - \frac {f(N)}{\log(N)}.
\end{equation}
Since $p$ is bounded from 0, then $S^*\leq O(\log(N))$
so (again using $N^{\frac {1}{\log(N)}}=2$)
$$cr=
\frac {c\log(S^*)}{\log(N)} \leq \frac {c\log \log(N)
+ O(1)}{\log(N)}.$$
All this is valid in particular when $1 \ll f(N) \ll
\log \log(N)$, so equation (\ref{eq__2}) becomes
\ifPersonal
$$ \frac 1 {k+1} \geq cr + \frac {f(N)}{\log(N)} =
\opo \frac{c\log \log(N)}{\log(N)}.$$ The last
condition is met for some $o(1)$ and any $k$ s.t.
$$k\leq \omo \frac {\log(N)}{c\log \log(N)}.\qedm$$ \else
$$ \frac 1 {k+1} \geq cr + \frac {f(N)}{\log(N)} =
\opo \frac{c\log \log(N)}{\log(N)}.\qedm$$ \fi

%=======================================================================
%         Upper Bound - Independence Number
%=======================================================================
%
The following upper-bound for the independence-number is larger than the bound of observation \ref{kwise_gives_exact_indp_num},
yet holds for significantly smaller densities $p$.

\btn {Independence-number upper bound}
\label{kwise_gives_up_bound_indp_set}
There exist constants $c_1,c_2$ \st for any \kwig \engnk
the following \as holds. There are no independent-sets
of size $S$ whenever either:
\begin{enumerate}
\item
$S \gg p^{-1/2}N^{3/4},~~k\geq 4$~~and~~$p \gg
N^{-1/2}$;~~or
\item
$S \geq c_1 \sqrt{\frac N p},~~k\geq \log(N)$~~and~~$p
\geq \frac {c_2 \log(N)}{N}$.
\end{enumerate}
\etn
%----------------------------------------------------
%
\noindent{\bf Proof.} By Theorem
\ref{kwise_gives_jumbl}, $\alpha$-jumbledness is \as
achieved. For item 1 we have $\alpha \gg
\sqrt{p}N^{3/4}$. For item 2 we have
$\alpha=O(\sqrt{pN})$. Then, any vertex set $U$
satisfies $e(U) \geq p|U|^2 - \alpha|U|$, so if $U$ is
independent, then $|U| \leq \alpha/p$. \qed

\ifPersonal Note that, in the context of the
independence number upper-bound, the power of the
jumbledness approach stands out more boldly compared
with the poor upper-bounds\footnote{Indeed, the
probability of having independent sets of size $S$ is
bounded by $\binom N S \Delta$, $\Delta$ being the \kwi
upper-bound on the probability of occurrence of a single
predetermined independent set.
Now, as $p \rightarrow 0$, then $\ln(1/1-p)\rightarrow
p$. Thus, when $p< (\log(N))^{-\Omega(1)}$, we have
$S^*= \Omega(1/p)$ and $\binom N S > (N/S)^{S} =
2^{(\log(N))^{\Omega(1)}}$, which is too large for
$\Delta$ to cancel when $k\leq \poly(\log(N))$.} that
the direct Chebyshev argument provides (thus, $k$-wise
independent only ensures the appropriate degrees,
co-degrees and number of 4-walks and then $I$ is
bounded very indirectly via jumbledness arguments).
\fi

%=======================================================================
%            $K$-Wise into Jumbledness Proof
%=======================================================================

\mysection{Appendix - $\bf k$-wise independence
guarantees optimal jumbledness}
\label{proof_kwise_gives_jumbl}

This appendix is dedicated to proving Theorem
\ref{kwise_gives_jumbl}.
Given an $N$-vertex graph $G$, consider the complete
graph $\bar G$, with
weight $1-p$ on any edge that appears in $G$, and
weight $-p$ on any other edge and on any self loop.
Let $A=A(\bar G)$ denote the corresponding $N \times N$ matrix
where $A_{u,w}=1-p$ whenever $u,w$ are adjacent in $G$ and
$A_{u,w}=-p$ otherwise (including the case $u=w$).

Let $\lambda=\lambda(\bar G)$ denote the largest eigenvalue
in absolute value
of $A$.
By the argument in \cite{ac} $G$ is $\lambda$-jumbled.
Indeed, for any two sets of vertices $U$ and $W$,
if we let $x_U$ and $x_W$ denote the characteristic vectors
of $U$ and $W$, respectively, then $x_U^t A x_W=e(U,W)-p|U||W|$
and as the $\ell_2$-norm of $Ax_W$ is at most $\lambda \sqrt {|W|}$,
and that of $x_U$ is $\sqrt {|U|}$,
it follows by Cauchy-Schwarz, that
$|e(U,W)-p|U||W||=|x_U^t A x_W| \leq \lambda \sqrt {|U||W|}$,
as needed.

Let $\Gamma= (v_0 \rightarrow v_1 \rightarrow \cdots
\rightarrow v_R= v_0)$ be an arbitrary closed walk
with $R$ steps in $\bar G$. Throughout, $\Gamma$ may
repeat vertices and edges and may traverse self-loops.
Let $W(\Gamma)=\prod_{j=0}^{R-1} A_{v_j,v_{j+1}}$, and
let $X=\sum_{\Gamma} W(\Gamma)$.

By Wigner's trace argument
\cite{wig}, for any graph distribution, and any even
$R\geq4$ and $\omega\gg 1$, \as
\begin{equation} \label{jumbl_vs_walks}
\lambda \leq (\omega \mathbb E (X))^{1/R}
\end{equation}
($\mathbb E(\cdot)$ stands for expectation). Thus,
establishing the desired jumbledness reduces to bounding
$\mathbb E= \mathbb E(X)$.
We first fix $t,j$, and bound the contribution to
$\mathbb E$ of a single walk, $\Gamma$, that traverses
exactly $t$ vertices and $j$ edges; later we bound the
number of walks with such $t,j$.
Let $\{e_1,...,e_{j}\}$ denote the set of all edges (excluding self-loops) used by $\Gamma$, where  $e_i$ is traversed precisely $q_i \geq 1$ times (we don't care how many times $e_i=\{u,w\}$ is traversed specifically from $u$ to $w$ or from $w$ to $u$).
As long as $k\geq R$ then the contribution of $\Gamma$ to $\mathbb E$ is bounded by $E (\Gamma)= \prod_{i=1}^{j}
[p(1-p)^{q_i} + (1-p)(-p)^{q_i}]$.
The latter equals $0$ if some $q_i=1$, so we focus on
walks where each $e_i$ is traversed at least twice.
Then,
\begin{equation} \label{exp_single_walk}
E(\Gamma) \leq \prod_{i=1}^{j} [p(1-p)^2 + (1-p)(-p)^2]
< p^{j}.
\end{equation}

{\bf Proving Theorem \ref{kwise_gives_jumbl}, item 1.}
Let $k=R=4$. There are only 2
types of walks: Walks with 3 vertices contribute
$O(p^2N^3)$ to $\mathbb E$, and walks with 2 (or 1)
vertices contribute only $O(pN^2)$, which is dominated
by the 3-vertex walks' contribution. By
(\ref{jumbl_vs_walks}), for any $\omega \gg 1$ \as
$\lambda \leq (\omega(p^2N^3))^{1/4}$. \qed~(Item 1)

{\bf Proving Theorem \ref{kwise_gives_jumbl}, item 2.}
We adopt the approach of F\"uredi-Komlos-Vu
\cite{fk2,vu}, who bound the number of walks with
given $(t,j)$, by encoding the walks in a 1:1 manner,
and then bound the number of code-words. We first
describe their encoding scheme (Section \ref{FKV}) and
later refine it (Section \ref{refined_encode}).

%--------------------------------------------------
\mysubsection{The F\"uredi-Komlos-Vu encoding}\label{FKV}
%--------------------------------------------------
%
Fix $\Gamma$ and consider the spanning-tree $T$ of
$\Gamma$, which consists of all the vertices visited by
$\Gamma$ and exactly those edges through which $\Gamma$
visits a vertex for the first time.
Edges (and consequently, steps) in $\Gamma$ are either
{\bf internal} ($e \in T$), or {\bf external} ($e \notin
T$). A step leading to a new vertex is called {\bf
positive}. A step traversing an internal edge for the
2'nd time is called {\bf negative}. Any other step is
called {\bf neutral}
(thus, all (+) steps are internal, and neutral steps are
either external, or pass through some internal edge for
the $i$'th time $i \geq 3$. Steps on self-loops are external).
The encoding of $\Gamma$ is composed of:
\begin{itemize}
\item
A list of all $t$ vertices visited by $\Gamma$, ordered
by their first appearance.
\item
A string of length $R$, where the $i$'th position
encodes the $i$'th step as follows. \subitem $\bullet$
Each positive step is encoded by (+). \subitem $\bullet$
Each negative step is encoded by (-). \subitem $\bullet$
Each neutral step $(u \rightarrow v)$ is encoded by
($v$).
\end{itemize}

How is $\Gamma$ retrieved from its encoding? The
starting vertex is known, since the order in which the
vertices appear in $\Gamma$ is known. Assuming that the
current position in the walk is known, then the next
position is also known if the next step is either
neutral or positive.
Ambiguity is possible only when we are about to traverse
a (-) step, and in addition the walk is currently at a
{\bf critical} vertex $x$. This means that the number of
internal edges $e_1,...,e_d$ that touch $x$, and have
been traversed exactly once (up to this point) is $\geq
2$.
For example, consider a walk starting with $1
\rightarrow 2 \rightarrow 3 \rightarrow 1 \rightarrow 4
\rightarrow 5 \rightarrow 1$. At this point, $x=1$ is
critical since both edges $e_1=\{1,2\}$ and
$e_2=\{1,4\}$ were traversed exactly once. If a (-) step
immediately follows, it is not clear to which $e_j$ this
current (-) refers. The encoding in \cite{fk2,vu} is
modified in some way \st critical steps can be decoded
un-ambiguously and the entire encoding-scheme becomes
1:1, as desired.
By Theorem 1.5 in \cite{vu} this suffices for proving
Item 2 in our theorem whenever $\Omega(\frac
{\log^4(N)} N) \leq p \leq 1 - \Omega(\frac
{\log^4(N)} N)$. \qed

For smaller $p$, we must refine the
original encoding significantly.

%-----------------------------------
\mysubsection{Our refined encoding}
%-----------------------------------
\label{refined_encode}
We start with a simple observation. Throughout the
analysis we set $k=R=\log(N)$, and let $\ell$ count
the number of external edges in $\Gamma$.
Let $\Phi_{t},\Phi^{\ell},\Phi_{t}^{\ell},$ resp. denote
the contribution to $\mathbb E$ of all walks with
exactly $t$ vertices, or exactly $\ell$ external edges,
or exactly $t$ vertices and $\ell$ external edges.
Clearly, $\mathbb E = \sum_{t} \Phi_{t}$, $\mathbb E =
\sum_{\ell} \Phi^{\ell}$, and $\mathbb E = \sum_{t,\ell}
\Phi_{t}^{\ell}$.
Since any of these sums has $\poly(\log(N))$ summands,
and since $R=\log(N)$ then $(\mathbb E)^{1/R}=
(1+o(1))\Phi^{1/R}$ whenever $\Phi$ bounds the maximal
term among all $\Phi_{t},\Phi^{\ell},\Phi_{t}^{\ell}$.
It therefore suffices to show that $\Phi \leq
(O(\sqrt{pN}))^R$.

We now give a high level description of our improved
analysis. We keep 3 ingredients from \cite{fk2,vu}:
(i) The entire (ordered) list of vertices is provided.
This contributes a $\Theta(N^t)$ multiplicative term to
the bound on $\mathbb E$.
(ii) Specific steps are encoded by a symbol from a fixed
alphabet (the original alphabet is
$\{+,-,\mathrm{neutral}\}$; our final alphabet will be
slightly larger). This contributes a $(\Theta(1))^R$
multiplicative term to $\mathbb E$.
(iii) Since $\Gamma$ traverses at least $t-1$ edges,
equation (\ref{exp_single_walk}) bounds the contribution
of each walk to $\mathbb E$ by $p^{t-1}$.
Thus, the combined contribution of (i)(ii) and (iii) to
the bound on $\lambda$ becomes (after taking the
$R$'th-root)
\begin{equation} \label{bound_basic_code}
\Theta((pN)^{t/R}).
\end{equation}

The latter partial encoding of (i)+(ii) is not 1:1, because
of the neutral and the critical steps.
Recall that there are $(t-1)$ edges in $T$, and since
(as mentioned above) all edges are traversed at least
twice, then $\Gamma$ has exactly $(t-1)$ (+) steps,
exactly $(t-1)$ (-) steps, and $m=R-2(t-1)$ neutral
steps.
In \cite{fk2,vu}, the trivial $t^m$ bound is used for
the contribution of the neutral steps to the total
number of code-words.
We strengthen \cite{fk2,vu}, mainly by the following
observations. First, we note that whenever $\ell$ is
`very large', then the entire contribution of all such
paths to $\lambda$ is negligible.
Next (and perhaps most significantly), we show that
whenever $\ell$ is not `very large', the following
holds: (i) Half of the neutral steps can be encoded very
economically, reducing the $t^m$ term from \cite{fk2,vu}
into roughly $t^{0.5m}$. (ii) All critical steps can be
encoded so economically, that their entire contribution
is (almost) dominated by that of the neutral steps.
Consequently, as $t \leq O(pN)$, we conclude that
$\lambda$ is bounded by (roughly) $\Theta(pN)^{\frac t
R} \times O(pN)^{\frac m {2R}}$ (the first term stems
from equation (\ref{bound_basic_code}), the second from
the neutral steps). Since $\frac m {2R} \sim 0.5 - \frac
t R$, the latter bound becomes $O(\sqrt{pN})$, as
desired. Details follow.

%----------------------------------------------------
\paragraph{Handling the `non-typical' walks (large $\ell$).}
%----------------------------------------------------
%
Clearly, the contribution of all $(t,\ell)$-walks to
$\mathbb E$ is bounded by $B= p^{\ell}(pN)^{t-1} t^R N$.
Indeed, ${N^t} t^R$ clearly bounds the number of walks,
and by equation (\ref{exp_single_walk}) the contribution
of each walk is bounded by $p^{(t-1)+\ell}$.
Let $\ell \geq 4\log\log(N)$. Since $t-1\leq 0.5 R,
R=\log(N),$ and $pN \leq O(\log^{4}(N))$, we have
$(pN)^{t-1} \leq \log(N)^{(4+o(1))(0.5\log(N))}=
\log(N)^{(2+o(1))\log(N)}$, $N=\log(N)^{o(\log(N))}$,
$t^R < \log(N)^{\log(N)}$, and $p^\ell \leq
(N^{1-o(1)})^{4\log\log(N)}=
\log(N)^{(-4+o(1))\log(N)}$.
Consequently $B\leq \log(N)^{(-4+2+1+o(1))\log(N)}\ll
1$. This concludes the treatment of the non-typical
walks.

%----------------------------------------------------
\paragraph{Handling the `typical' walks (small $\ell$).}
%----------------------------------------------------
%
A new encoding is required to handle the typical walks.
As before, the encoding includes the names of all $t$
vertices, ordered by their first appearance, and all (+)
and all non-critical (-) steps are simply encoded by (+)
and (-) and decoded trivially.
Our new perspective is thinking of the entire walk as
composed of sequences of internal steps separated by
external steps. We first encode the external steps
economically, and prove that this enables to
economically encode the critical steps as well.
Next, we handle the internal sequences. We first provide
some general observations regarding arbitrary internal
sequences. Then, we describe how to handle the specific
case of encoding a closed internal sequence. Finally, we
generalize the latter to encoding open internal
sequences as well.

%----------------------------------------------------
\mysubsubsection{Encoding external steps.}
%----------------------------------------------------
%
\label{encode_external_step}
To exploit the small number of external edges, we add
the following to the code.
\begin{itemize}
\item
A list of all $\ell$ external edges
${e_1,...,e_{\ell}}$.
\item
Each external step on $e_i$ is encoded by (${i,d}$),
where the bit $d$ specifies the direction in which $e_i$
is traversed.
\end{itemize}
Thus, encoding a singe external step has only $2\ell=
\Theta(\log\log(N))$ possible values. This improves
upon \cite{fk2,vu} where such steps are encoded by
their end-vertex which might have $\Theta(\log(N))$
possible values.

%----------------------------------------------------
\mysubsubsection{Encoding critical steps.}
%----------------------------------------------------
%
\label{encode_critical_step}
\ifPersonal Our approach
is different then \cite{fk2,vu}. While they
(expensively) encode an entire sequence of (-) steps
that starts at a critical step, we (economically)
encode only the critical step itself.
\fi
Recall that a step $\bar s$ is critical in $\Gamma$, if
$\bar s$ is taken from a vertex $x$, \st that $x$ has $d
\geq 2$ {\bf critical edges} for $\bar s$. Critical
edges are internal edges $e_1,...,e_d$ that touch $x$
and have been traversed exactly once up to $\bar s$.
We will show that each $e_i$ can be associated with a
unique external edge $e$.
\ifPersonal
(we do not claim that all critical edges in $\Gamma$ can
be associated uniquely with an external edge. But the
current critical edges at step $\bar s$ can).
\fi
This will enable us to encode $\bar s$ using $e$ which
has only $\Theta(\log\log(N))$ possible values.
\ifPersonal
This eventually improves upon \cite{fk2,vu}, where any
sequence of (-) steps that starts at a critical step
$\bar s$ is encoded by the name of some vertex which
might have $\Theta(\log(N))$ possible values.
\fi Specifically,
consider any critical edge $e_i=\{x,w\}$ for $\bar s$
which is not the first edge leading to $x$ in $\Gamma$.
Consider the step $s_i$ where $e_i$ is traversed for the
first time in $\Gamma$. Since $x$ had already appeared
in $\Gamma$, then $s_i= (x \into w)$.
If we omit $e_i$ from $T$, we partition $T$ into 2
disjoint sub-trees: $T_1$ which contains $x$ and $T_2$
which contains $w$. Since $e_i$ is critical then the
first time we return to $x$ (after $s_i$), is not via
$(w \into x)$.
Thus, there must exist some external edge $e$ (that
connects $T_2$ to $T_1$) that is traversed between $s_i$
and the first time we return to $x$. We call the first
of these edges $e$ the {\bf external criticality edge
(ECE)} of $e_i$, and denote it by $c(e_i,\bar s)$.
Clearly, different $e_i$-s have distinct ECEs, and in
addition, at step $\bar s$ all ECEs are well defined by
previous steps in $\Gamma$. Consequently, the following
encoding is un-ambiguous.
\begin{itemize}
\item
Let $\bar s$ be a critical step from $x$, with critical
edges $e_1,...,e_d$, and external criticality edges
$c(e_i,\bar s)$. If $\bar s$ traverses the first edge
that leads to $x$ in $\Gamma$ we encode $\bar s$ by (-).
Otherwise, we encode $\bar s$ by the position of
$c(e_i,\bar s)$ in the list of external edges.
\end{itemize}
Note that not-critical negative steps are encoded by (-)
as before. This concludes the treatment of external and
critical steps. We finally encode neutral internal
steps.

%----------------------------------------------------
\mysubsubsection{Encoding an arbitrary sequence of
internal steps.}
%----------------------------------------------------
%
Let $S=(s_1,...,s_q)$ be a `maximal' sequence of
internal steps. Here maximality means that (i) either
the step previous to $s_1$ was external or that $s_1$ is
the first step in the entire walk, and that (ii) either
the next step after $s_q$  is external, or that $s_q$ is
the last step in the entire walk (maximality does not
mean that there are no longer internal sequences $S'$ in
$\Gamma$).
We remark that, in general, some of the edges used by
$S$ may have been traversed before $S$ started and some
may be introduced by $S$ for the first time.

Let $x$ be the starting vertex of $S$. Fix some vertex
$w \neq x$ visited by $S$ and let $u=u(w)$ denote the
{\bf predecessor} of $w$ in $S$ (so the first time we
reach $w$ in $S$ is via $(u \into w)$).
Clearly, after each time we step $u \into w$, then the
only way to return to $u$ is by stepping $w \into u$
(otherwise we get a cycle from $u$ to itself in $T$).
Thus, when we pass $e$ for the $j$'th time during $S$ we
go forward $(u \into w)$ when $j$ is odd and go backward
$(w \into u)$ when $j$ is even. We call this the {\bf
forward-backward observation}.
Since the predecessor is uniquely determined by previous
steps in $\Gamma$, we can encode backward steps very
economically.
\begin{itemize}
\item
A neutral-backward step is (economically) encoded by
($nb$).
\item
A neutral-forward step ($u \into w$) is (explicitly)
encoded by ($nf,w$).
\end{itemize}
Given this, we desire to demonstrate that many of the
neutral steps in $S$ go backward. We first handle the
following simple case.

%----------------------------------------------------
\myParagraph{8.2.3.1~~~Encoding a closed sequence.}
%----------------------------------------------------
Assume $S$ is closed, namely, the end vertex of $s_q$
is the starting vertex, $x$, of $s_1$. We claim that
at least half the neutral steps in $S$ go backward. We
actually prove the latter for every edge $e$ in $S$.
Let $\# f(e)$ and $\# b(e)$ resp.~denote the number of
forward and backward steps on an arbitrary edge $e$
during $S$. Note that currently, not only neutral but
also (+) and (-) steps are counted. We show that $\#
f(e)=\# b(e)$.
Indeed, otherwise, by the forward-backward observation
the last step on $e$ was a forward step $s_i=(u \into
w)$, and clearly there exists a path from $x$ to $u$ in
$T$. However, since $S$ is closed there must exist
another path in $T$ from $w$ to $x$ that avoids stepping
$(w \into u)$ - a contradiction.
Now, let $\# nf(e)$ and $\# nb(e)$ count the number of
neutral-forward and neutral-backward steps on $e$ during
$S$. By the above, there are at least 2 steps on $e$.
There are 3 cases:
(i) If $e$ was never used prior to $S$ the first step is
$(+f)$, the second is $(-b)$. The next steps (if any
exist) come in pairs of $(nf)(nb)$. (ii) If $e$ was used
at least twice prior to $S$, then all steps come in
pairs of $(nf)(nb)$. (iii) If $e$ was traversed exactly
once prior to $S$, the first 2 steps are $(-f)(nb)$, and
all consequent steps (if any exist) come in pairs of
$(nf)(nb)$. In cases (i),(ii) we get $\# nf(e)=\# nb(e)$
and in case (iii) $\# nf(e)=\# nb(e)-1$.
Anyway, at least half of the neutral steps in a
closed-internal sequence can be encoded economically.
This concludes our analysis for closed sequences.

The problem is that for open sequences, $S$, it might
hold that all $m(S)$ neutral steps in $S$ are forward,
and by the \cite{fk2,vu} encoding-scheme these steps
contribute a huge $t^{m(S)}$ factor to the bound on the
number of code words. To overcome this, we use the
following.

%----------------------------------------------------
\myParagraph{8.2.3.2~~~Encoding an open sequence.}
%----------------------------------------------------
Let $x \neq y$ denote the start-vertex and end-vertex of
some open maximal internal sequence $S$. Clearly, in $T$
there exists a unique path $P=(x=x_1 \into x_2 \into ...
\into x_r=y)$.
All the steps in $S$ can be uniquely partitioned into
2 categories. (i) Steps that traverse $P$ (either
forward ($x_i \into x_{i+1}$) or backward ($x_{i+1}
\into x_{i}$)). (ii) Entire sub-sequence $S'$, where
each $S'$ starts and ends at some path-vertex $x_i$,
but never touch $P$ at any other vertex other than
$x_i$. Such $S'$ is a closed internal sequence and is
encoded as discussed in Section $8.2.3.1$.
We first modify the encoding simply \st path steps are
explicitly encoded as such.
\begin{itemize}
\item
Each positive-path step is encoded by ($+p$).
\item
Each forward-neutral-path step is encoded by ($npf$).
\item
Each backward-neutral-path step is encoded by ($npb$).
\item
Each forward-negative-path step is encoded by ($-pf$)~
(if the step $\bar s$ is critical and has an external
criticality edge, the index of this edge in the list
of external edges is added to the encoding of $\bar s$
as in Section \ref{encode_critical_step}).
\item
Each backward-negative-path step is encoded by ($-pb$).
\end{itemize}

Clearly, {\bf if} the entire path $P$ is known, then the
latter encoding suffices to decode any path-step,
because on a path there is a unique forward-step and
unique backward-step from each vertex. The question is
how to recover the path itself.
First, the end vertex $y$ of $P$ is well defined by the
encoding. Indeed, if $S$ is immediately followed by an
external step $\bar s$, than the encoding of $\bar s$
determines $y$. Otherwise, $y$ is the last vertex in
$\Gamma$, which is the (already known) first vertex of
$\Gamma$.
To recover the remaining vertices in $P$, we call a
vertex on $P$ either {\bf old} or {\bf new} according to
whether it appeared in $\Gamma$ prior to $S$ or not. If
all vertices are old, since we know $x$ and $y$, and
since there is a unique path connecting $x$ to $y$ in
$T$, then $P$ is uniquely defined.
Otherwise, some vertices in $P$ are new. We claim that
no new vertex is followed in $P$ by an old vertex.
Otherwise, the path $P$ includes a step $x_i \into
x_{i+1}$ where $x_{i}$ is new but $x_{i+1}$ is old. This
means that before $S$ started there was a path in $T$
from $x_{i+1}$ to $x_1$ (at each point in the walk,
there exists a unique sub-tree of $T$ that spans all
the vertices traversed so far). Thus, $P$ closes a cycle
in $T$ from $x_{i+1}$ to itself - a contradiction.
Therefore, if there are any new vertices there exists a
unique {\bf final old} vertex $\bar x$ along $P$. If
$\bar x$ is known, then the path from $x$ to $\bar x$ is
unique (since there is a unique path between any
vertex-pair in $T$). In this case, the other part of $P$
from $\bar x$ to $y$ is also well defined, because it
consists only of new vertices (recall that the order in
which new vertices appear in $\Gamma$ is explicitly
encoded). This covers all possible cases.
Note that actually, if all steps on $P$ are (+) and (-),
then they are already uniquely decodable as before.
Thus, the only addition required for decoding path steps
is:
\begin{itemize}
\item
Let $S$ be an open internal sequence, with path $P$
that contains at least a single new vertex and at
least a single neutral step. Let $\bar x$ be the final
old vertex in $P$. Then the symbol $\bar x$ is added
to the encoding of the first neutral path step in $S$.
\end{itemize}
The main benefit here is that instead of encoding the
end-point of each forward neutral step, it suffices to
encode once the entire `direction' of the path (this
approach is similar to the \cite{fk2,vu} encoding of
critical steps).

%----------------------------------------------------
\myParagraph{Wrapping up.}
%----------------------------------------------------
By all the above, the final encoding (including all
aforementioned modifications) is 1:1 as desired. We
currently fix any $\ell < 4\log\log(N)$ and bound the
contribution $E_{t,\ell}$ of all $(t,\ell)$-walks to
$\mathbb E$. Specifically, we bound the contribution
of various parts in our encoding to $E_{t,\ell}$. Each
contribution introduces a new multiplicative term to
$E_{t,\ell}$.
As before,
\begin{itemize}
\item
Choosing the (list of ordered) $t$ vertices to appear in
$\Gamma$ contributes $(N)_t=\Theta(N^t)$.
\item
The basic encoding of each step as some combination of
positive/negative/neutral path/non-path forward/backward
contributes $(\Theta(1))^R$.
\item
The contribution of each single walk $\Gamma$ is
$\Theta(p^{t-1+\ell})$.
\end{itemize}
We now consider the critical and neutral steps. Recall
$m$ is the total number of neutral steps. Let $m_1$
count the number of external steps. Let $m_2$ count the
neutral steps in closed internal sequences. Let $m_3$
count the number of open internal sequences $S$ \st
their path $P=P(S)$ contains at least a single neutral
step and at least a single new vertex.

\begin{itemize}
\item
Choosing the ${\ell}$ external edges to appear in
$\Gamma$ contributes $\binom {t^2}{\ell}=
O(t^{2\ell})$.
\item
By Section \ref{encode_external_step} encoding the
external steps contributes ${(2\ell)}^{m_1}$.
\item
By Section \ref{encode_critical_step} encoding the
critical steps contributes at most $(\ell+1)^{t-1}$.
\item
By Section 8.2.3.1 encoding the neutral steps in closed
internal sequences contributes at most $t^{0.5 m_2}$.
\item
By Section 8.2.3.2 encoding the neutral steps on the paths
of open internal sequences contributes at most
$t^{m_3}$.
\end{itemize}
Recall that to bound $\lambda$ we are about to take the
$R$'th root of $\mathbb E$, and that we are willing to
tolerate small $\Theta$ factors in the bound on
$\lambda$.
Since there are only $(\log(N))^{\Theta(1)}$ possible
$t,\ell,m_1,m_2,m_3$, and since
$(\log(N))^{\Theta(\frac 1 R)}=1+o(1)$, then we may
consider only the choice of $t,\ell,m_1,m_2,m_3$ that
maximizes the bound (on the contribution to $\mathbb
E$).
In addition, we may (i) Ignore the $(O(1))^R$ factor
from encoding specific steps as a combination of
$\{+,-,n,p,f,b\}$. (ii) Consider $N^{t-1}$ instead of
$N^t$ (because $N^{\frac 1 R}=2$). (iii) Ignore the
$t^{2\ell}$ factor from the choice of external edges
(since $t^{2\ell} < \log(N)^{8\log\log(N)}=
2^{o(R)}$). (iv) Replace the $(\ell+1)^{t-1}$ term
with $(\ell)^{t}$.

Combining all the remaining (un-ignored) terms yields
the following expression
$$\Psi= (pN)^{t-1} t^{0.5{m_2}+m_3} p^{\ell}
\ell^{m_1+t}.$$
As ${m_1+t} < {\log(N)}$, and $p \geq N^{-1+o(1)}$,
then $p^{\ell} \ell^{t+m_1}= 2^{-(1-o(1))\ell\log(N)}
2^{+\log(\ell)\log(N)} \leq 1$.
Next, consider any open internal sequence $S$ counted in
$m_3$. For any such $S$ (except possibly the last one),
there exists a unique neutral-external step that
terminates $S$, so $m_3 \leq 0.5(m-m_1)+1$. Thus,
$0.5{m_2} + m_3 \leq 0.5m+1= 0.5(R-2(t-1))+1$. Thus,
$t<pN$ implies $\Psi< (pN)^{t-1} t^{0.5(R-2(t-1))+1} <
(pN)^{0.5R+1}.$ All the above gives $\lambda \leq
\Theta(1)\Psi^{1/R}= \Theta(\sqrt{pN})$. \qed

\end{document}